\documentclass[11pt]{article}
\usepackage{amsmath,amssymb,bbm,pslatex,xcolor,graphics}

\usepackage{xcolor}
\definecolor{blue}{rgb}{0,0,1}
\usepackage[pagebackref,bookmarks]{hyperref} 
\hypersetup{pagebackref=true,bookmarks=true,
       hyperfigures=true,
	unicode=true,
	colorlinks=true,
	citecolor=blue,
	linkcolor=blue,
	anchorcolor=red
}
\usepackage[textsize=footnotesize,color=blue!40, bordercolor=white]{todonotes}

\usepackage{etoolbox}
\makeatletter
\def\new@mathgroup{\alloc@8\mathgroup\mathchardef\@cclvi}
\patchcmd{\document@select@group}{\sixt@@n}{\@cclvi}{}{}
\patchcmd{\select@group}{\sixt@@n}{\@cclvi}{}{}
\makeatother

\usepackage{amsmath}
\usepackage{amssymb}

\usepackage{graphicx,bbm,color,calc,ifthen,theorem,capt-of,tocloft}

\usepackage{fancyhdr}

\newcommand{\call}{ {\mathcal{L}}}

\newcommand{\phiabp}[2]{\ensuremath{\Phi_{\beta +1}^{(  \alpha )}}}

\newcommand{\nod}{{\noindent}}

\newcommand{\ul}{\ensuremath{\ulcorner}}
\newcommand{\ur}{\ensuremath{\urcorner}}

\newcommand{\bu}{\ensuremath{\bullet}}
\newcommand{\nat}{\ensuremath{\mathbbm{N}}}
\newcommand{\re}{\ensuremath{\mathbbm{R}}}

\newcommand{\cant}{2\textsuperscript{\ensuremath{\mathbbm{N}}}}

\newcommand{\emp}{{\varnothing}}
\newcommand{\pa}[1]{\ensuremath{\langle #1 \rangle}}

\newcommand{\vp}{\ensuremath{\varphi}}

\newcommand{\da}{{\downarrow}}
\newcommand{\ua}{{\uparrow}}
\newcommand{\la}{{\langle}}
\newcommand{\ra}{{\rangle}}
\newcommand{\dom}{\text{dom}}
\newcommand{\back}{{\backslash}}
\newcommand{\power}{{\mathcal{P}}}

\newcommand{\ie}{{\itshape{i.e.}}{\hspace{0.25em}}}
\newcommand{\eg}{{\itshape{e.g.}}{\hspace{0.25em}}}
\newcommand{\etc}{{\itshape{etc.}}{\hspace{0.25em}}}
\newcommand{\via}{{\itshape{via }}{\hspace{0.25em}}}
\newcommand{\cf}{{\itshape{cf. }}{\hspace{0.25em}}}
\newcommand{\mm}{{\itshape{m.m.}}{\hspace{0.25em}}}

\newcommand{\all}{{\forall}}
\newcommand{\Equi}{{\,\Longleftrightarrow\,}}
\newcommand{\equi}{{\,\longleftrightarrow\,}}
\newcommand{\ex}{{\exists}}

\newcommand{\sset}{ { \subseteq } }
\newcommand{\lr}{{\leftrightarrow}}
\newcommand{\imp}{{\,\longrightarrow\,}}

\newcommand{\rem}{{\noindent}{\bfseries{Remark: }}}
\newcommand{\fin}[1]{\ensuremath{[  ]^{< \omega}}}

\newcommand{\pf}{{\noindent}{\textbf{Proof: }}}
\newcommand{\dfs}{\ensuremath{=_{\ensuremath{\operatorname{df}}}}}

\newcommand{\oddpagetext}[1]{\newcommand{\pageoddheader}{{\small }}}
\newcommand{\evenpagetext}[1]{\newcommand{\pageevenheader}{{\small }}}

\def\sset{\mbox{ $\subseteq$ }}
\catcode`\<=\active \def<{
\fontencoding{T1}\selectfont\symbol{60}\fontencoding{\encodingdefault}}
\catcode`\>=\active \def>{
\fontencoding{T1}\selectfont\symbol{62}\fontencoding{\encodingdefault}}
\catcode`\|=\active \def|{
\fontencoding{T1}\selectfont\symbol{124}\fontencoding{\encodingdefault}}

\newcommand{\tmstrong}[1]{\textbf{#1}}

\newcommand{\nosymbol}{}
\newcommand{\nobracket}{}
\newcommand{\nocomma}{}

\newcommand{\tmem}[1]{{\em #1\/}}
\newcommand{\tmop}[1]{\ensuremath{\operatorname{#1}}}

\newcommand{\tmtextit}[1]{{\it{#1}}}
\newcommand{\tmverbatim}[1]{{\ttfamily{#1}}}
\newcommand{\tmtextbf}[1]{{\bfseries{#1}}}
\newtheorem{theorem}{Theorem}[section]

\newtheorem{corollary}[theorem]{Corollary}
\newtheorem{definition}[theorem]{Definition}

\newtheorem{lemma}[theorem]{Lemma}
\newtheorem{remark}[theorem]{Remark}

\def\nod{\noindent}

\def\mb{\mbox{}}

\def\nod{\noindent}

\def\ed{\end{document}}

\def\bth{\begin{theorem}}
\def\eth{\end{theorem}}

\def\blem{\begin{lemma}}
\def\elem{\end{lemma}}

\def\bdefn{\begin{definition}}
\def\edefn{\end{definition}}

\newcommand{\qed}{\mbox{ } \hfill Q.E.D.}

\def\mid{\, | \, }

\def\m``{\mbox{``}}

\def\le0{L[E_{0}]}

\def\ii{\iota}

\def\mm{{\em mutatis mutandis}}

\def\power{{P}} 

\def\egc{\em{e.g.},}

\def\ran{\mathop{ran}}
\def\ul{\raisebox{-0.1em}{$\ulcorner\!$}}
\def\ur{\raisebox{-0.1em}{$\!\urcorner$}}

\def\ul{{\ulcorner\!}}
\def\ur{{\!\urcorner}}
\def\widebar{\bar}
\usepackage{times}

\begin{document}

\date{17.iii.20}

\title{Characterisations of Variant Transfinite Computational Models: Infinite Time Turing, Ordinal Time Turing, and Blum-Shub-Smale machines}
\author{P.D.~Welch}
\maketitle
\begin{abstract}
  We consider how changes in transfinite machine architecture can sometimes
  alter substantially their capabilities. We approach the subject by answering
  three open problems touching on: firstly differing halting time
  considerations for machines with multiple as opposed to single heads,
  secondly space requirements, and lastly limit rules. We: 1)
  use admissibility theory, $\Sigma_{2}$-codes and $\Pi_{3}$-reflection properties in the constructible hierarchy to
   classify the halting times of
  ITTMs with multiple independent heads; the same for Ordinal Turing Machines
  which have $\tmop{On}$ length tapes; 2) determine which admissible lengths
  of tapes for transfinite time machines with long tapes allow the machine to
  address each of their cells - a question raised by B. Rin; 3) characterise
  exactly the strength and behaviour of transfinitely acting Blum-Shub-Smale
  machines using a $\tmop{Liminf}$ rule on their registers - thereby
  establishing there is a universal such machine. This is in contradistinction to 
  the machine using a `continuity' rule which fails to be universal.
\end{abstract}

\section{{\small{}}Introduction}

The study of transfinite computational models fully reemerged with the
appearance of the Infinite Time Turing machines (ITTMs) of Hamkins and Kidder
in {\cite{HL}}\footnote{We should like to thank Merlin Carl, Lorenzo Galeotti
and Philipp Schlicht for conversations concerning some of the issues and
arguments here.}\footnote{Keywords: 
Turing Machine, computability, Blum-Shub-Smale, generalized recursion, admissible set}. 
It was not that there were no precursors to the idea of
transcending the finite by some forms of recursive model: Machover
{\cite{Mach61}}, and Takeuti {\cite{Tak60}} formulated equational calculi
involving ordinal recursive set functions, as did more famously perhaps Jensen
and Karp {\cite{JeKa71}}, and Gandy {\cite{Ga74}}. It is hard not to view
equational calculi of this type as being simply abstract `machines' as well.
However Platek explicitly thought of a register-like machine containing
transfinite, rather than finite, ordinals. Kleene himself used a model
comprising infinite but wellfounded trees whose nodes had (ordinary) Turing
machines attached, for use in computing higher type recursion
{\cite{Kl62a}},{\cite{Kl62b}}. \ Rogers in his fundamental text {\cite{Rog}}
used the analogy of the `$\aleph_{0}$-mind' as having some countably infinite
memory storage with an ability to sort and manipulate countable amounts of
information in single steps. This was in his account of a description of the
hyperarithmetic sets. Even ITTMs had a precursor (it later transpired) with a
device invented by Burd {\cite{Burd84}}, for which, although being entirely
similar to ITTMs, he seemingly only wrote one routine which took just $\omega
\cdot 2$ steps to complete.

However {\cite{HL}} caught the imagination of not a few people, and thereafter
various authors examined all sorts of finite based machines and tried to
propel them into the transfinite. Others sought to tweak the parameters of an
ITTM. Restricting to the latter for the moment, it transpired that the model
was very robust under structural changes: whether one considered
$\liminf$rather than the original $\limsup$ rule at limit stages made no
difference; neither did relaxing the requirement that the read/write (R/W)\footnote{We use a variety of initialisms: ``R/W'' for ``read/write''; ``OT'' for ``output tape''; ``TM'' is not ``Turing machine'' but ``Theory Machine''.}
head return to the origin and enter a special limit state $q_{L}$ make any
difference to the functions on Cantor space, {\cant}, so computed. We further
considered the question as to whether some other limit rule, say another
$\Sigma_{2}$-definable limit rule, would allow for a wider class of functions
to be rendered ITTM-computable. In {\cite{W}} we showed that essentially the
Liminf rule was `maximal' or `universal' amongst all other $\Sigma_{2}$-rules.
This with hindsight is not very surprising: once one has established that
standard ITTMs can compute (codes for) the levels of $L$ up to $L_{\zeta}$
where $\zeta$ is the least $\Sigma_{2}${\tmem{-extendible ordinal}}, that is
the least so that there is some $\Sigma > \zeta$ with $L_{\zeta}
\prec_{\Sigma_{2}} L_{\Sigma}$ (see the ``$\lambda$-$\zeta$-$\Sigma$-Theorem
of {\cite{W}}), then, with any other limit rule being absolutely definable to
$L$, this shows that an ITTM equipped with a different rule can not do better
than computing codes up to $L_{\zeta}$, and would be subject to the same
looping behaviour at $[ \zeta , \Sigma ]$ as the standard ITTM. It is then
easy to see that the standard ITTM can simulate, {\via}the $L$-hierarchy if
need be, the non-standard one. So the latter cannot be stronger. (Limit rules
at higher degrees of definability, so $\Sigma_{3}$ and above, are altogether a
different matter.)

Hamkins and Seabold ({\cite{HaSe}}) considered whether the three Input, Output
and Scratch tapes of an ITTM could be replaced by a single tape.  It turned out that the class of computable functions $F:2^{\omega}
\imp \omega$ was unchanged but $F:2^{\omega} \imp 2^{\omega}$ would be
smaller. Indeed, while it is easy to easy
enough to conceive of using recursive subclasses of a single tape for scratch
and output,  it was shown that the compression of a
tape to a single `output tape' for presentation of the resulting element of
$2^{\omega}$ as output would go awry. It was noted that with larger alphabets,
or a special reserved cell or some other gadget, would restore the full class.
An example of this can be found in {\cite{W2010}}  where we considered the idea of using cells with blank contents:
this allowed a blank for ``ambiguity'' to be used at limit stages $\lambda$
when the cell's value had changed unboundedly often below $\lambda$. \ \ On a
one-tape model this restored the full computational power of the three-tape
ITTM. Certain considerations would still nevertheless mean that halting times,
or lengths of gaps in halting times would still change if the tapes or
alphabet were to change - as one might imagine. \

The original formulation of ITTM's in \cite{HL} used the idea that at limit ordinal times, $\mu$ say, the machine entered a special
`reserved state' $q_{L}$ to indicate it was at a limit, but also the R/W head was always returned to the starting position
so that it was reading the first cells of each of the three tapes. Subsequent discussions showed that a designated state $q_{L}$ was somewhat unnecessary, as was the requirement to return to the beginning: instead the R/W head can be returned to the $\liminf$ of the head positions prior to the limit
ordinal $\mu$, and the state number $i$ (of $q_{i}$) at that time $\mu$ was then set to be the $\liminf$ of the previous state numbers. This seemed anyway more rational as with reasonable programming desiderata, this arrangement had the R/W head, and the program, entering the chief or head subroutine that was initiated unboundedly often in time $\mu$. It is easy to check that, \eg\!\!, the class of ITTM computable functions, semi-decidable sets, \etc\!\!, is unchanged by this. The arguments of \cite{HL} that showed admissible ordinals cannot be halting times of an ITTM program, and that only admissibles may start gaps in these halting times (see Def. \ref{23} below)  of \cite{W09} adapt to this situation.
We have adopted this formulation throughout this paper.

We address the question here of ITTMs with not only multiple tapes, but
multiple and independently moving R/W heads. \ This of course has to be
programmed into the state or transition table of the Turing program, as each
head will be moving according to the independent diktats of whatever the head
is observing at the time. However this is quite unproblematic and familiar to
Turing machine theorists. A complete snapshot of the machine at any stage of
time $\alpha$ now must also give information about the position of each of the
heads, in addition to the $\omega$-sequence of cell values on each of the
tapes. To be concrete we imagine below an ITTM with two independent \ R/W
heads and at least two tapes. This suffices for the arguments below: more
tapes and heads does not make a difference to the results. The results
illustrate a difference with the single head models considered to date: the
same classes of functions are thereby computed, but the spectrum of halting
times, whilst retaining the same supremum as before, now has a different
distribution. For example $\omega_{1}^{\tmop{ck}}$ now becomes a halting time
of a multihead ITTM.

These results are also planned with the {\tmem{Ordinal Turing Machines}} of
Koepke {\cite{K05}} in mind. These OTMs were a natural extension, in that a
tape was allowed with an $\tmop{On}$-length of cells. Koepke gave a standard
definition of such a machine using a single tape. However then he fairly
quickly went to also allowing multiple tape, multi-head models. All such
computed (on zero input) the same class of (codes of) sets: namely all of $L$
the G\"odel constructible hierarchy. Thus a new presentation of this hierarchy
became available which only used this hardware plus a standard Turing
program. Indeed he phrased this in terms of computing a truth predicate for
$L$, where a formula and ordinal parameters could be submitted {\via}the input
tape, for query as to its truth in $L$ and 1/0 output would be calculated.
Questions of `halting times' for programs $P_{e}$ for this model were not
considered by him.

\begin{definition}
  \label{23} We say that an ordinal $\alpha$ is $\tmop{OTM}${\tmem{-clockable}}
  {\tmem{(ITTM-clockable)}} if for some $\tmop{OTM}$-program ({\tmem{ITTM}}
  program) $P_{e} $, some $k \in \omega$, we have $P_{e} ( k )
  \da^{\alpha}$ , that is, it halts in exactly $\alpha$ steps, and thus
  $\alpha$ is the {\tmem{halting time}} of $P_{e} ( k )$.
\end{definition}

Carl in {\cite{Carl20}} considers this for a multitape and multihead version of
OTM. He shows: any $\Sigma_{2}$-admissible ordinal is not clockable in this
sense. This provides an upper bound. As a lower bound he shows that any `gap'
in the clockable ordinals (meaning a maximal interval $[ \gamma , \delta )$
where no program $P_{e} ( k )$ on integer input $k$ halts at any time $\alpha
\in [ \gamma , \delta )$) must be initiated by an ordinal $\gamma$ which is at
least an admissible \ limit of admissibles. The possibility remained open that
an admissible limit of admissibles could be such a gap-starting $\gamma$. We
show that this is not the case, and eliminate the space between upper and
lower bounds by pinning down the exact gap-starting ordinals by showing in
Section 2:\\

\nod{\tmstrong{{\nod}Theorem \ref{Th2.3}}} {\tmem{ Suppose we consider multiple head {\tmem{ITTM}}s.

(a) Let $\alpha$ be $\Pi_{3}${\tmem{}}-reflecting. Then $\alpha$ is not
$\tmop{ITTM}$-clockable.

(b) The least $\Pi_{3}${\tmem{}}-reflecting $\alpha_{0}$ starts the first gap
in such clockables. In general, if $\alpha$ is not
$\Pi_{3}${\tmem{}}-reflecting, then $\alpha$ does not start a gap; 
 hence for such $\alpha$, if it is a limit of clockables, it is itself clockable.

 The same holds true for {\tmem{OTM}}s if again we allow multiple
heads.}}\\

{\nod}{\tmem{{\tmstrong{}}{\tmstrong{{\tmem{Corollary{\tmstrong{}}}}{\tmem{
\ref{C2.8}}}}} }}{\tmem{If $\gamma$ is a supremum of halting times in either
the \tmem{ITTM} or \tmem{OTM} multihead models,}}{\tmem{ then}}  $\gamma$ {\tmem{starts a gap iff
$\gamma$ is $\Pi_{3}$-reflecting.}}\\

Hence in contradistinction to the single head ITTMs, $\omega_{1}^{\tmop{ck}}$
is a halting time for the multihead version. The above then also settles then the
question of the clockable OTM ordinals, and those starting gaps. Nothing so
far answers the following:\\

\nod {\bf{Question}}: {\tmem{Is $\omega_{1}^{\tmop{ck}}$ clockable by a
single-head \tmem{OTM} program?}}\\

{\nod}However we'll see as a Corollary of Theorem 2.10 that $\omega_{1}^{\tmop{ck}}$ is also a halting time
of a single-head OTM.\\

{\nod}{\tmstrong{Corollary \ref{C2.11}}} {\tmem{Let $\gamma$ be a supremum of
single-head {\tmem{OTM}} halting times. Then:}}

\ \ \ \ \ \ \ \ \ \ \ \ \ \ \ \ $  \gamma$ {\tmem{starts a gap $\Equi$
$\gamma$ is either $\Pi_{3}$-reflecting, or is an admissible limit of
$\Sigma_{2}$-extendibles}}.\\

{\nod}Here we recall:

\begin{definition}
  An ordinal $\xi$ is a $( \Sigma_{2} )${\tmem{-extendible ordinal}}, if
  $L_{\xi}$ has a proper transitive $\Sigma_{2}$-end-extension to some
  ${L}_{\sigma}$: $L_{\xi} \prec_{\Sigma_{2}} L_{\sigma}$.
\end{definition}

The $\tmop{OTM}$s with ordinal tapes, as the ITTMs even do, beg the question
of the nature of ITTM versions with tapes of some set length $\alpha$, with
$\omega < \alpha$. It is reasonable to consider $\alpha$ with some amount of
closure, for example, under primitive recursive set or ordinal functions, or
else to be admissible. In the early days of ITTM or study of ``$\alpha$-ITTM''
study, there had been some hope that since a machine theoretic description of
these levels of the $L$ hierarchy was possible, and the actions of such
machines were so slowly restricted to one-step-at-a-time actions, that
possible some new properties of $L$ would discovered, or some approach to a
new, or even the older, Jensen-style fine structure might emerge. Perhaps even
a new proof, or aspects of proof, of combinatorial principles such as Jensen's
Square principle might be possible? Silver machines {\cite{Ri79}},
{\cite{De}}, {\cite{BL80}} give one a hope in this direction. However it
seemed that these hopes were to be disappointed. Jensen fine structure (and
Silver machines) depend very much on the notion of $\Sigma_{1}$-definability.
\ The $\Sigma_{2}$-limit rule of $\alpha$-ITTMs works against any finer
dissection of the Jensen fine structure, or proofs of $\Box$. Arguments
constructing $\Box$-sequences depend on taking $\Sigma_{0}$-hulls and unions
of such. A $\Sigma_{2}$-hull is simply too rich.\\

There is a broader context in which to discuss ITTMs and their variants, and
that is of the notion of {\tmem{quasi-inductive definition}}. 
This is strictly speaking not a discussion needed to undertake the results of the paper, but is more motivational of some the viewpoints we hold concerning basic definitions, and gives some historical background as well as tying in the older, but very developed, theory of inductive definitions, with ITTM theory.
This notion is taken
from nomenclature devised by John Burgess {\cite{B}} in a paper discussing the
Herzberger variant of the {\tmem{revision theory of truth}} of Gupta and
Belnap {\cite{GB}}. This turns out to be formally equivalent, in some loose
sense, when performed as a truth theory over the natural structure $\nat =
\pa{\nat ,0,' ,+, \times}$ and the theory of ITTMs. 

This revision theory uses a recursive process to define truth sets of sentences, so ``truth-sets'' over a first order structure, such as $\nat =
\pa{\nat ,0,' ,+, \times, T}$ in, here, the language of arithmetic augmented by a $\dot T$-symbol  for a one place predicate whose extension $T_{\alpha+1}$ at stage $\alpha+1$ is the set of G\"odel numbers of sentences in this extended language true in $\pa{\nat ,0,' ,+, \times, T_{\alpha}}$. Unlike the Kripkean theories of truth that build up partial truth sets in such a language monotonically, here the extension is supposed to be {\em total}. However liar sentences (such as a simple liar $L\equi \neg L$) and other self-referential sentences are possible as elements of $T_{\alpha}$ - and the rules for negation require that $\ul
L\ur \in
T_{\alpha}\equi \ul L \ur \notin T_{\alpha +1}$; hence the process is non-monotonic. The question that exercised revision theorists was what to do a limit stages, thus how should \eg $T_{\omega}$ be defined? Herzberger took a $\liminf$ rule: $T_{\lambda}= \liminf_{\alpha\rightarrow \lambda}T_{\alpha}$ and the reader can no doubt surmise some connections to ITTMs building up the extension of the cell values of the Turing machine over stages $\alpha$. Often for a Herzberger revision sequence that incorporated this liminf rule,
one would start with $T_{0}=\emp$ although the more sophisticated revision theories  considered all possible starting distributions of truth sets, and indeed all possible limit rules for a stage $\lambda$, with the proviso that it be consistent with the set of sentences already stabilized below $\lambda$.
 
Burgess considered
{\tmem{arithmetical quasi-inductive definitions}} which were in general
intermediate between the idea of recursive quasi-inductive definition afforded
by an ITTM, and the fully first order, or $\Delta_{\omega}^{0}$-, idea
inherent in a Herzberger revision sequence. Just as at limit stages an ITTM
defines its cell values by a $\Sigma_{2}$-limit liminf rule, so does a
Herzberger sequence define extensions to the $\lambda$-th stage $T_{\lambda}$
of a truth predicate, by a Liminf rule of the previous `revised' sets of
G\"odel numbers $T_{\alpha}$ for $\alpha < \lambda$. L\"owe first pointed out
the similarities between ITTMs and these revision theories. In {\cite{L}} he
wrote an ITTM program for simulating a given Herzberger revision sequence
where the length of that sequence was itself writable by an independent ITTM
program.  That this latter restriction was unnecessary was argued in
{\cite{W1}}: ITTMs can fully simulate any such revision sequences without
fixing a length in advance.

Such $\liminf$ rules do not have many examples in either mathematical, or
philosophical, logic. This is perhaps unsurprising. The inbuilt nature of the
ordinals to their definition, and the fact that the ordinals required turned
out to need $\Pi^{1}_{3}$-$\tmop{CA}_{0}$ to define them is a rarity (\cf \cite{Si99}). This
contrasts with the clearly natural {\tmem{inductive definition}} in particular
the {\tmem{positive elementary inductive definition}} of $\cite{M}$.

A set $S \sset | \mathfrak{A} | \nobracket \nobracket$ (for some arbitrary
relational structure $\mathfrak{A}$) is {\tmem{inductive}} over  $\mathfrak{A}$ in this
sense if it is (1-1) reducible to a fixed point arising from a positive
elementary induction definition; more particularly it is a {\tmem{section}} of
that fixed point.

In some more detail:
let $\mathfrak{A} =  \pa{A,R_{1} , \ldots \, ,R_{k}}$ be a relational
structure. The language
 $\call_{\mathfrak{A}}  $ contains constant symbols
$\dot{p}$ for elements $p$ of $A= | \mathfrak{A} | $, and relation symbols $\dot{R_{1}} , \ldots \nocomma
\dot{R_{}}_{k}$. Let $\varphi \in \call_{\mathfrak{A}}$ be $\varphi ( v_{0} ,
\ldots v_{n-1} ,p_{1} , \ldots p_{m} ,S )$ be a formula in which $S$ occurs
{\tmem{positively}}. We let $\varphi$ define an {\tmem{operator}}:\\
$$\Gamma_{\!\varphi}  : \power ( A^{n} ) \imp \power ( A^{n} )\mbox{ by
} \Gamma_{\!\varphi} ( S ) = \{ \bar{x} \,: \varphi (
\bar{x} , \bar{p} ,S )^{\mathfrak{A}} \} .$$

\nod Such an inductive definition by the positivity of the occurrence of $S$ results in $\Gamma_{\vp}$ being a monotone operator, and hence we are guaranteed a fixed point.
We thus define iterates:

$$I^{< \xi}_{\varphi} = \bigcup_{\tau < \xi} I^{\tau}_{\varphi}\, \quad \, ; \quad I^{\xi}_{\varphi} = \Gamma _{\vp}( I^{< \xi}_{\varphi} )
\quad ; \quad I_{\varphi}
\dfs I^{\infty}_{\varphi}.$$

\begin{definition}
  (i) $P$ is {\em{inductive over $\mathfrak{A}$}}, 
 if for  some $\vp$, for some $\bar{p} \in A^{n}$, then $$P= \{ \bar{x}  
  \mid   ( \bar{x} , \bar{p} ) \in I_{\vp} \} \hspace{2em}\mbox{ ``}P \in
  \tmop{IND}_{\mathfrak{A}}\mbox{''}$$
   (ii) $P$ is {\em{hyperelementary over $\mathfrak{A}$}} if  both
  $P$ and $\neg P$ are inductive. \quad ``$P \in \tmop{HYP}_{\mathfrak{A}}$''.
\end{definition}

\nod The inductive sets are closed under: $=,R_{i} , \vee , \wedge ,
\ex^{\mathfrak{A}} , \all^{\mathfrak{A}}$, $\tmop{HYP}$- substitution.
We are going to replace the $\Sigma_{1}$-rule of simple unions at limit stages
by the $\Sigma_{2}$-rule of taking $\liminf$'s:
$$I^{< \xi}_{\varphi} = 
\liminf_{\alpha < \xi} I^{\alpha}_{\vp} 
= \bigcup_{\alpha < \xi} \bigcap_{\alpha
< \tau < \xi} I^{\tau}_{\vp}\quad ; \quad I_{\varphi} = I^{\infty}_{\varphi}
\dfs \liminf_{\tau < \infty} I^{\tau}_{\vp} .$$

\begin{definition}
  (i) $P$ is {\em{quasi-inductive over $\mathfrak{A}$}},
   if for some $\vp$, for some $\bar{p} \in A^{n}$, then $$P= \{ \bar{x}  
  \mid   ( \bar{x} , \bar{p} ) \in I_{\vp} \} \hspace{2em}\mbox{ ``}P \in
  \tmop{qIND}_{\mathfrak{A}}\mbox{''}.$$
   (ii) $P$ is {\em{quasi-hyperelementary over $\mathfrak{A}$}} if  both
  $P$ and $\neg P$ are quasi-inductive.\quad ``$P \in \tmop{qHYP}_{\mathfrak{A}}$''.
\end{definition}

\nod The quasi-inductive sets are then also closed under: $=,R_{i} , \vee , \wedge ,
\ex^{\mathfrak{A}} , \all^{\mathfrak{A}}$, $\tmop{qHYP}$- substitution.


That book (\cite{M}) is dedicated to investigating the mathematics of inductive definitions
as they occur over quite arbitrary abstract structures (or at least those
equipped with a minimal amount of pairing and coding machinery - the
`acceptable' structures). This defines a class of subsets of the domain of the
structure, as well as the class of sets of such subsets, as the `inductively
definable' ones. Theorems such as Stage Comparison, Boundedness Lemmata, an
Abstract Spector-Gandy Theorem {\tmem{etc., etc.}} all come into place.

It is surprising how much of this can be extended when one deploys the
notion of elementary quasi-inductive definition as defined above.
The parallels with ITTMs are clear: consideration of the
tapes contents moves us to think of a computation as {\tmem{convergent}} if
the contents of the output tape {\tmem{eventually settles}} to some infinite
sequence of values. Thinking of the ITTM with its tapes and cell structure as
itself the abstract structure $\mathfrak{A}$ over which we perform the
quasi-induction, a set of integers is {\tmem{quasi-inductive}} if it is again
a{\tmem{ section}} of a set arising as a quasi-induction in this sense over \
$\mathfrak{A}$ (this means in the case of an $\omega$-ITTM, that the set is
(1-1) reducible to a fixed part of the tapes' contents at stage $\zeta$). It
is {\tmem{quasi-hyperelementary}} if it is both quasi-inductive and co-quasi-inductive.
This comes down to: the quasi-hyperelementary sets are those that are eventually
settled (on the OT), {\ie} those that are {\tmem{eventually writable}} also known as
{\tmem{(eventually) decidable}}. \ The quasi-inductive sets are the
{\tmem{(eventually) semi-decidable}} sets.

Thus the emphasis of the ITTM or $\alpha$-ITTM, is shifted away from the
halting computations where the computation formally halts, to those that have
simply finished writing to their OT's - although they may be fiddling away
pointlessly on their scratch tapes: in short formally halting is just a
special case of {\tmem{convergence}}. This reflects the analysis of the
$\lambda$-$\zeta$-$\Sigma$-Theorem: the prime ordinal of interest is $\zeta$,
not $\lambda$: the latter is merely the supremum of the halting times of
programs on integer input; $\zeta$ is the supremum of the times when a
program on input may have ceased writing to its OT. Without analysing this
phenomenon of ceasing to write, and $\zeta$, the proof of the
$\lambda$-$\zeta$-$\Sigma$-Theorem would not have occurred. We thus regard
`ceasing to write' as prior to the notion of halting for ITTMs. An
over-emphasis on `halting' as opposed to this eventual settled behaviour,
perhaps lead the original authors into investigating a degree theory in
{\cite{HL}} that was based on analogies with ordinary Turing degree theory
with its halting problems. This, to the current author's thinking, was a red
herring: the analogy of ITTM-degrees is not with Turing degrees but to
something nearer hyperdegrees, indeed to something intermediate between
hyperdegrees and $\Delta^{1}_{2}$, say $\Delta^{1}_{3 \back 2}$-degrees. The
mathematics of ITTM-degrees bears out that this is the correct analogy: the
semi-decidable sets form a {\tmem{Spector pointclas}}s; there is a
{\tmem{Spector criterion}} for the degrees and their jump operation. The above
notion of convergence of ITTMs fits exactly with the notion of quasi-inductive
over a Turing machine (or over $\nat$). Indeed both Herzberger revision
sequences and ITTMs provide the prime and first examples of quasi-inductive
processes over $\nat$.\\

Pursuing this we come to the second question addressed here. When putting the
framework of positive elementary inductive definitions in place for an
arbitrary abstract structure $\mathfrak{A}$, Moschovakis worked in a language
$\call_{\mathfrak{A}}$ containing constants $c_{a}$ for every $a$ in the
domain $| \mathfrak{A} | \nobracket \nobracket$. One clearly has to be able to
refer to the elements of the structure in order to build up extensions of some
predicate: which are in, which are out.

B. Rin {\cite{Rin}} identifies a problem in the model of Turing machines with
transfinite length tapes being unable to access all of their cells if one
disbars `constants', that is parameters, from their architecture. This is
indeed a problem. Parameters are normally allowed as part of the definition of
a transfinite computation, particularly in OTMs: we allow a computation to
proceed from a finite number of $1$'s as input in cells $C_{\xi_{1}} \, ,
\ldots \, \nocomma \, ,C_{\xi_{n}}  $ to indicate, or stand in for, a
finite list of ordinals $\xi_{1} , \ldots \nocomma , \xi_{n}$. Such parameters
then function as constants allowing us to address, or as Rin defines, `reach'
each cell. We may `address' a cell on an $\alpha$-machine if we have a name
for that cell ({\eg} in the case of $\alpha = \omega$) say, or we may `reach'
a cell by having a particular program $P_{e} ( 0 )$ that halts exactly on that
cell. \ \ (In which case case we might call `$e$' a name for that cell and say
that this cell is `addressable' or `reachable'.) For $\alpha = \omega$ we
clearly need no such parameters to name or reach cells. \ But even for
countable $\alpha$ the problem arises of what happens when $\alpha$ is
sufficiently big that the addressing of all the cells on the tape is
impossible? If $\alpha = \omega_{1}$ then computing on such a machine would
allow the computation of, say a function $F: \omega_{1} \imp 2$ on the work
tape, but one would only have access to countably many values: namely those $F
( \beta )$ for $\beta  $ addressable in this way, if one wanted to work with
the function.

The main question left unresolved in {\cite{Rin}} is the identity of the
following $\delta_{0}$:\\

{\nod}{\bf Question} {\tmem{What is the least countable $\delta_{0}$ so that there is some
cell $C_{\alpha}$ for an $\alpha < \delta_{0}$, which cannot be reached by a
program run on a $\delta_{0}$-ITTM machine?}}\\

We characterise below precisely the lengths of such tapes for which this
problem occurs: they are those for which there is no computation during which
a wellordering appears which collapses the length of the tape, here $\delta$,
to be countable.\\

{\nod}{\tmstrong{Theorem }}{\tmem{{\tmem{{\tmstrong{\ref{Th3.1}{\tmem{}}}}}} \
\ Every cell $\alpha < \delta  $ can be reached by a $\delta$-ITTM if and only
if during a run of some program on such a machine, a wellorder of $\omega$
appears on its tape at some stage which has order type $\delta$.}}\\

\nod{\tmstrong{{\nod}Corollary}} {\tmem{{\tmstrong{{\tmem{\ref{Cor3.4}}}}}}}
{\tmem{\tmverbatim{}Every cell can be reached by a $\delta$-ITTM \ if and only
if $L_{\Sigma ( \delta )} \models$``$\delta  $ is countable'' (where $\Sigma (
\delta )$ is the supremum of ordinals  coded as subsets of $\omega$  appearing on the tapes of
some $\delta$-ITTM on integer input, at some stage of its computation).}}\\

In {\cite{CaRiSc18}} it is established that $\delta_{0}$ is the least ordinal
which is uncountable in $L_{\Sigma ( \delta_{0} )}$, answering the question
above. However the proof here is simpler and direct, and moreover
characterises cell reachability for all primitive recursively closed ordinals $\delta$.

We move next in Section 4 to Infinite Time Blum-Shub-Smale machines
(IBSSM's). The author recalls a conversation in New York with Hamkins and
Koepke, as to how to go about propelling BSS-machines into the transfinite.
Two limit rules were proposed: the ``continuity'' rule, that required one to
have a continuous limit of each register's contents at a limit time; this
seemed the most promising. A ``Liminf'' rule was mooted, but there was overall
no further discussion. Koepke and Seyfferth then established several facts in
{\cite{KoSe12}} concerning the `Continuity'-IBSSM's, the principal one being
that any such machine on integer/rational input would either halt or be in a
permanent loop by time $\omega^{\omega}$. Thus the continuity rule was a
stringent one. This left open the strength of the functions computable by such
machines. In {\cite{We2015}} (Thm. 11) we showed that the class of functions
so computed was exactly those in $L_{\omega^{\omega}}$. In fact it showed
more, that there were equivalent formulations in terms of classes of
functions, with ``polynomial time'' ITTM-computable functions (these were
those that halted in time some ``polynomial in $\omega$'', or in short by some
time $\omega^{k}$ for a $k< \omega$), and those functions on $\omega$-strings
generated by the ``safe recursive set functions'' of Buss et al.
{\cite{BeBuFr12}}. We shall not define this class here, but refer the reader
instead to {\cite{BeBuFr12}}. \ One thus has some confluence of models
computing the same classes - a sub-Church's thesis perhaps for polytime on
$\omega$-strings (or equivalently elements of Cantor space $\cant$).\\

{\tmem{{\tmem{{\tmstrong{Theorem }}({\cite{We2015}}}} Thm. 11) The following
classes of functions of the form $F: (2^{\mathbbm{N}} )^{k} \rightarrow
2^{\mathbbm{N}}$ are extensionally equivalent:

(I) Those functions computed by a continuous {\tmem{IBSSM}} machine;

(II) Those functions that are polynomial time {\tmem{ITTM:}} those computed by
some time $\omega^{k}$ for a $k< \omega$;

(III) Those functions that are safe recursive set functions.}}\\

Section 4 here seeks to show that the strength of IBSSM's equipped with a
Liminf rule for the register contents, rather than the original continuity
requirement allows, unsurprisingly, for much greater computational power. It
turns out that they have the same power as ITTMs no less. Perhaps this is
more surprising given that the notion of Liminf in a Euclidean setting of
$\re$ with its usual metric, is rather different from a Liminf of cell values
in the setting of ITTMs on Cantor space $\cant$. It shows that Liminf rules
subsume much of the successor step character of the processes involved.

Essentially this comes from the fact that any ITTM computation $P_{e} ( 0 )$
can be simulated on a \ `Liminf'-IBSSM $B_{\widebar{e}} ( 0.0 )$. Consequently the
IBSSM's can have no lesser strength than ITTMs. Since ITTMs can compute
(codes for) levels of the $L_{\alpha}$ hierarchy this means IBSSM's can do the
same. It is easier to see that any IBSSM computations can be considered as
absolute to $L$. This implies that they can be simulated by ITTMs. Hence we have the
same $\lambda$-$\zeta$-$\Sigma$ phenomenon for IBSSM's, as for ITTMs: namely
they must loop at latest by time $\zeta$ with a loop repeating at time
$\Sigma$. Another corollary (Cor. \ref{Cor20}) is that Liminf-IBSSM's have a
universality property - this is because ITTMs do and the two classes of
machine are `bi-simulable'. There is thus a universal such Liminf-IBSSM.
This contrasts with the `Continuity'-IBSSMs where there is no such universal machine.
(This follows easily from Theorem 1 of \cite{KoSe12} which shows that for any IBSSM machine using the continuity criterion at
limit stages, if a program has $k$ computation nodes in its flow chart, then any halting computation must have length less than 
$\omega^{{k+1}}$, and so is not reaching all the functions IBSSM-computable. Applying that, it is immediate that no such program can be universal.)\\

{\nod}{\tmstrong{Corollary}} \ref{C4.5} {{\em For
{\tmem{$\tmop{liminf}$-IBSSM}}'s the set $Z$ of reals on which an \tmem{IBSSM}
computation (on rational input) is convergent, is precisely the set $Z $ of
those reals in $L_{\zeta}$.}}\\

(A {\tmem{convergent computation}} (on rational input) is one in which a
designated register, {\eg}$R_{0}$, has a final settled value from some point
onwards.)\\

\nod{\em \bf Preliminaries}\\

We assume the reader is familiar with the machine architectures under
discussion - to limit space we do not review these here, but instead refer the
reader to the original papers or to \cite{Carl2019}. 

The results here very much tie in with low level set theory and G\"odel's constructible hierarchy
of sets $L= \bigcup_{\alpha\in On}L_{\alpha}$ the basic construction of which we assume (v. \cite{D}
for example or elsewhere).

The theory of ITTM halting times, and indeed computable sets is intimately tied up with the levels of the hierarchy which are {\em admissible models} of {\em Kripke-Platek set theory} which we may take as the weakening of \tmem{ZFC} by restricting to $\Delta_{0}$-Collection and $\Delta_{1}$-Separation. It is well known that in this theory $\Sigma_{1}$-Replacement is provable.
The reader can consult
{\cite{Bar}} for these results, or for admissibility theory in general. We shall need in addition to this notions of {\em reflecting} levels of the hierarchy.

\begin{definition}
  An ordinal $\alpha > \omega$ is $\Pi_{n}${\tmem{-reflecting}}, if for any
  $\varphi ( v_{0} ) \in \Pi_{n}$, for any $x \in L_{\alpha}$, if $\varphi [ x
  ]^{L_{\alpha}}$ then for some $\beta < \alpha$ $\varphi [ x ]^{L_{\beta}}$.
\end{definition}

{\rem} (1) Admissibles are $\Pi_{2}$-reflecting, and hence it is easily seen
they are reflecting for Boolean combinations of $\Pi_{2}$ and $\Sigma_{2}$
formulae. Clearly even with parameters allowed, if a formula reflects once because
of such a principle holding at $\alpha$, it will reflect unboundedly often
below $\alpha$.

(2) A $\Pi_{3}${\tmem{}}-reflecting ordinal is admissible, and is a limit of
such. It is easy to see that (a) $\Sigma_{2}$-admissibles are
$\Pi_{3}${\tmem{}}-reflecting; however (b) the first admissible limit of
admissibles is not; (c) for any $n<\omega$, the first $\Pi_{n}$-reflecting ordinal is less than the first
$\Sigma_{2}$ admissible.

(3) The axioms of $\tmop{KP}$ are formulable in a $\Pi_{3}$-manner. Hence we
can always take the $\beta$ in the definition above to be admissible when
considering $\Pi_{3}$-reflection. (A $\Pi_{3}$-reflecting ordinal is
`recursively Mahlo' in the terminology of {\cite{PH2}} - although not
conversely.)\\

 At one or two points we refer to, and indeed use the ``Theory Machine'' (TM) (see \cite{FrWe08}). 
To make this paper a little more self-contained we sketch this program here. The idea of the TM is
to produce in a coherent, and sequential fashion, {\em codes} (meaning reals coding the $\in$-diagrams of levels of the $L$-hierarchy, $L_{\alpha}$, together with their $\Sigma_{2}$-{\em theories}.

The motivating point behind the Theory Machine is that it constructs such codes and theories
 in a smooth uniform fashion. In \cite{W} we used, in a slightly {\em ad hoc} fashion, a program that ran the universal ITTM, and simply collected together sums of ordinals appearing on the tapes as it progressed. Given such a sum, in a side man\oe uvre, we could then construct a segment of the $L$-hierarchy along this ordinal sum for our inspection. These sums, almost by definition, stretched out unboundedly in $\Sigma$. What the {TM} does is cut this collection process of ordinals out, and simply produces the $L_{\alpha}$ for $\alpha < \Sigma$ more directly.

An important point from the theory machine, is that in essence it is a version of the universal infinite time Turing machine program: ITTM operations are absolute to $L$, and the theories of various levels of $L_{\alpha}$ contain the information about the universal ITTM run up to stage $\alpha$. Then the information of the universal ITTM (on integer inputs) is all implicitly contained in the course of computation of the one single TM.\\

We set $T^{2}_{\alpha}$ to be the $\Sigma_{2}$-theory in the language of set theory of $(J_{\alpha},\in)$. We use
the $J$-hierarchy rather than the traditional $L$-hierarchy to avail ourselves of uniformly definable Skolem functions.
This is not terribly important, but using the $L$-hierarchy is a bit more
awkward. Recall that $On\cap J_{\alpha}= \omega\cdot \alpha$ and if $\omega \cdot \alpha = \alpha$ (certainly when $\alpha$ is admissible or primitively recursively closed) then $J_{\alpha}
=L_{\alpha}$ in any case, hence replacing $J_{\alpha}$'s by $L_{\alpha}$'s in the sequel will be no great lie for those unfamiliar with the former. The reader need not study the proof of the next Lemma nor the one following in order to understand the arguments of the paper.

\begin{lemma}\label{L1.3}
  There is an \tmem{ITTM} programme $P_{e} =$ {\em TM} which does not converge, but
  continuously produces alternately codes $x_{\alpha}  $ for levels
  $J_{\alpha}$ and their $\Sigma_{2}$-theories $T^{2}_{\alpha}$ for $\alpha <
  \Sigma$. At stage $\Sigma$ as $T^{2}_{\Sigma} =T^{2}_{\zeta}$ the {\em TM} loops
  back and reproduces the code $x_{\zeta}$ and continues this process
  thereafter repeating codes and theories for $\alpha \in [ \zeta , \Sigma )
  .$ 
\end{lemma}

{\pf} We briefly sketch the effective procedure to be formalised.
 The input to
$\tmop{TM}$ is presumed to be zero. 
On the output tape at round $\alpha$ a code for $x_{\alpha}$ is first written in one designated area and from this the theory $T^{2}_{\alpha}$ is calculated and written in another. Then $x_{\alpha+1}$ is calculated from  $T^{2}_{\alpha}$ and overwrites $x_{\alpha}$. Then $T^{2}_{\alpha+1}$ is determined, and overwrites $T^{2}_{\alpha}$. At a limit stage $\lambda$ the area for codes contains in general just an overwritten  mess of integers, but the area for the theories contains - by the usual Liminf rules $\widehat{T}_{\lambda}\dfs \liminf_{\alpha \rightarrow \lambda}T^{2}_{\lambda}$. If $\lambda$ is admissible then this equals $T^{2}_{\lambda}$. But even if not, one may show that, uniformly in $\lambda$, $T^{2}_{\lambda}$ is recursively enumerable in the liminf theory $\widehat{T}_{\lambda}$. From this we may construct a code for $J_{\lambda}$ and continue. Each round of this process only takes roughly speaking another $\omega^{2}$-many steps.  As one can surmise, the process must cycle around when it reaches stage $\Sigma$ and this is because we must have $\widehat{T}_{\zeta}= \widehat{T}_{\Sigma}$ and the next code produced falls back to $x_{\zeta}$ once more.\\ 

This description suffices for this paper, and, again, the rest of this proof can be ignored, but we give a few more details here for the interested reader.
 A theory $T$ is written with $\varphi_{n} \in
T$ iff the $n$'th cell contains a 1.  If $\sigma$ is a sentence in this language, we let $\ul \sigma\ur$ 
be that $n$ where $\sigma$ is $\vp_{n}$ in the given recursive enumeration. 
In the first $\omega^{2} + \omega \cdot
2$ stages $\tmop{TM}$ writes the code of $J_{1} =L_{\omega} = \tmop{HF}$ and \
its $\Delta_{0}$-diagram which we shall denote $d_1$,  to two reserved tapes, and its $\Sigma_{2}$-theory to
the OT. (It takes less than this, but it keeps the induction
bookkeeping straight.) We assume inductively that at time $\omega^{2} \cdot
\alpha + \omega \cdot 2$ the OT contains the theory $T^{2}_{\alpha}$ of
$\pa{J_{\alpha} , \in}$ and the reserved tapes again the $\Delta_{0}$-diagram
of $J_{\alpha}$, $d_{\alpha}$, and a code for $J_{\alpha}$. With the theory
$T_{\alpha}^{2}$ of $J_{\alpha}$ TM can construct a code $x_{\alpha +1}$ for
$J_{\alpha +1}$ in $\omega^{2}$ additional steps together with its
$\Delta_{0}$-diagram $d_{\alpha +1}$ (this is essentially an exercise, but uses the fact that for $\alpha < \Sigma$ there is a uniform parameter free $\Sigma_{2}^{J_{\alpha}}$-definable partial map $f_{\alpha}:\omega \rightarrow J_{\alpha}$ which is onto).

 We are are now at stage 
$\omega^{2} \cdot \alpha + \omega \cdot 2+ \omega^{2}$. \ In an additional
$\omega \cdot 2$ steps $T^{2}_{\alpha +1}$ is calculated from $d_{\alpha +1}$
and written to OT. This will take us to stage $\omega^{2} \cdot (
\alpha +1 ) + \omega \cdot 2$. (This all takes some routine work to make
clear, but essentially this two quantifier theory can be recovered from the double Turing jump of the code $d_\alpha$.
Each jump can be written out by an ITTM in $\omega$-steps
(in fact the double jump can be so written, but we can ignore that), thus
requiring $\omega \cdot 2$ steps to write out the two jumps and thus obtain
the complete theory $T^{2}_{\alpha+1} .$)

Of course we do this writing simply by changing the cells one by one according
to what has appeared or disappeared passing from $T^{2}_{\alpha}  $ to
$T_{\alpha +1}^{2}$. If $\varphi_{n}$ is in both theories, then the 1 in the
$n$'th cell is not changed to a 0 and then back again to a 1. \ By this method
of writing, at a limit stage $\omega^{2} \cdot \lambda$ for $\tmop{Lim} (
\lambda )$, $\widehat{T}_{\lambda}$ is on the OT, and thus the true
$T^{2}_{\lambda}$ is r.e. in the OT, by the next lemma. Hence in $\omega$
further steps it can then write the correct $T^{2}_{\lambda}$ to the OT, thus
by stage $\omega^{2} \cdot \lambda + \omega$.

The following is then Lemma 1 from \cite{FrWe08}.
\begin{lemma}\label{L1.4}
  There is an (ordinary) Turing recursive function $f: \omega \times \omega
  \imp \omega$, so given by an index $e$, so that for any $\tmop{Lim} (
  \lambda )$ \ satisfying $J_{\lambda} \models$``Every set $x$ is countable'',
  if we set $T= \widehat{T}_{\lambda}$ \ then $T^{2}_{\lambda}$ is uniformly
  r.e. in $T$, {\via}$f$, that is: \ \  $\ul \sigma \ur \in
  T^{2}_{\lambda} \equi \ex i f \left( i, \ul \sigma \ur \right) \in T$. \ \ 
\end{lemma}

{\pf} Let $\sigma \equiv \ex u \psi ( u )$ with $\psi \in \Pi_{1}$. Let $h_{\tau}$ be the canonical $\Sigma_{1}$ Skolem function for $J_{\tau}$. 
If $\pa{\vp_{i}\mid i<\omega}$ enumerates the $\Sigma_{1}$ 
formul\ae \,\!  of the language of set theory, then if $(\ex x \vp_{i}(x, v_{1},\,\ldots\, ,\, v_{n}))^{J_{\tau}}$ then 
$(\vp_{i}(h_{\tau}(i, \pa{v_{1},\,\ldots\, ,\, v_{n}}, v_{1},\,\ldots\, ,\, v_{n}))^{J_{\tau}}$.
Thus $h_{\tau}$ is a partial function from $\omega\times \mb^{<\omega}J_{\tau}$ to $J_{\tau}$
and has the same definition in over any $J_{\tau}$ uniformly in $\tau$. 
Let $S^{1}_{\tau}\dfs \{ \gamma < \tau \mid J_{\gamma}\prec_{\Sigma_{1}}J_{\tau}\}$. Then for any $\gamma\in S^{1}_{\tau}$ we have that $h_{\tau}$``$\omega \times \{
\gamma\}$ is the least $\Sigma_{1}$ Skolem Hull $X\prec_{\Sigma_{1}}J_{\tau}$ with $\gamma\in X$. If $J_{\tau} \models$``Every set $x$ is countable'', then it is easy to show that $\gamma +1\subset X$ and moreover $X\cap On\cap J_{\tau}$ is the least element of $S^{1}_{\tau}\cup \{\omega \cdot \tau\} $ strictly above $ \gamma$.\\

\nod {\tmem{Claim: $\sigma \in T^{2}_{\lambda} \equi \ex i \left[ \ex \tau_{0}
\all \tau \in ( \tau_{0} , \lambda ) \right.$}}
{\tmem{$ J_{\tau} \models$}}``$\ex \beta
\in S^{1}_{\tau} \left( ( \sigma )^{J_{\beta}} \vee \left( h_{\tau} ( i, \beta
) \da \wedge \psi [ h_{\tau} ( i, \beta ) ]^{J_{\tau}} \right)
\right)$''$\nobracket ].$\\


{\pf} of Claim.

{\em Case 1. $S^{1}_{\lambda}$ is unbounded in $\lambda$.}

\nod Suppose the left hand side holds of $\sigma$. Suppose $\psi ( u_{0}
)^{J_{\lambda}}$ holds for $u_{0}$. Then for some sufficiently large $\beta \in
S^{1}_{\lambda}$, $u_{0} \in J_{\beta}$, and then $\psi ( u_{0}
)^{J_{\beta}}$. But $\beta \in S^{1}_{\lambda} \imp \beta \in S^{1}_{\tau}$ \
for any $\lambda >\tau > \beta$; consequently the first disjunct of the right hand side
holds. For the converse direction, fix the given $i$. By the Case hypothesis
we can assume that $\tau$ itself is in $S^{1}_{\lambda}$. But then if the
first disjunct holds, if $( \sigma )^{J_{\beta}}$ and $\beta \in S^{1}_{\tau}$
then \tmtextbf{$\beta \in S^{1}_{\lambda}$} and thence $( \sigma
)^{J_{\lambda}}$. If the second disjunct holds for the supposed $i$ $\psi [
h_{\tau} ( i, \beta ) ]^{J_{\lambda}}$ holds for the same reasons.

{\em Case 2 $\beta_{0} \dfs \max  S^{1}_{\lambda} < \lambda$ exists.}

\nod By the bullet points above every $x \in J_{\lambda}$ is of the form
$h_{\lambda} ( i, \beta_{0} )$. Again suppose the left hand side holds of
$\sigma$ and $\psi ( u_{0} )^{J_{\lambda}}$ holds for $u_{0}$. In particular
now $u_{0} =h_{\lambda} ( i, \beta_{0} )$ for some $i$. Let $\tau_{0} \geq
\beta_{0}  $ be sufficiently large so that $\left( h_{\tau_{0}} ( i, \beta )
\da \right)^{J_{\tau_{0}}}$ and thence, by the fact of $\psi  $ being
$\Pi_{1}$, $( \psi [ h_{\tau_{0}} ( i, \beta ) ] )^{J_{\tau_{0}}}$. By the
upwards persistence of $\Sigma_{1}$ formulae in the first case and downwards
persistence of $\psi$ in the second case, these will hold in all larger
$J_{\tau}$ for $\tau \leq \lambda$ replacing $\tau_{0}$. But now the second
disjunct of the right hand side holds.

Conversely suppose the right hand side holds. Let $i$ be as supposed. By the
maximality of $\beta_{0}$ for unboundedly many $\tau' \in ( \beta_{0} ,
\lambda )$ some new $\Sigma_{1}$-sentence about $\beta_{0}$ becomes true first
in $J_{\tau' +1}$. Pick such a $\tau = \tau' +1$ of this form. Such a $\tau$
ensures that $S^{1}_{\tau} =S^{1}_{\lambda}$ and thence $\max  S^{1}_{\tau} =
\beta_{0}$ too. So suppose the first disjunct holds for such a successor
$\tau$. Then if $( \sigma )^{J_{\beta}}$ holds for a $\beta \in S^{1}_{\tau}
=S_{\lambda}^{1}$ we shall have $( \sigma )^{J_{\lambda}}$ and we are done.
Thus we now suppose the first disjunct fails for $\tau$ of this form; pick any
such $\tau$, then the second disjunct holds as witnessed by a $\beta \in
S^{1}_{\tau}$.

Then if $\beta < \beta_{0}$ then $\left( \ex y ( y=h_{\tau} ( i, \beta ) )
\right)^{J_{\tau}}$ implies $\left( \ex y ( y=h_{\beta_{0}} ( i, \beta ) )
\right)^{J_{\beta_{0}}}$ by the uniformity of the definition of the
$\Sigma_{1}$-skolem function $h$, and the fact of $\beta_{0} \in
S^{1}_{\tau}$. But then $\left( h_{\beta_{0}} ( i, \beta ) \da \wedge \psi [
h_{\beta_{0}} ( i, \beta ) ] \right)^{J_{\beta_{0}}}$. But this entails that
the first conjunct holds for $\tau$, which we are assuming does not happen.
Hence we must have $\beta = \beta_{0}$. \ However then $\psi [ h_{\tau} ( i,
\beta_{0} ) ]^{J_{\tau}}$ for any $\tau$ of this form, and so for such $\tau$
arbitrarily large below $\lambda$. By the upwards persistence of $h_{\gamma} (
i, \beta_{0} )$ for $\gamma \in [ \tau , \lambda ]$ we have \ $\psi [ h_{\tau}
( i, \beta_{0} ) ]^{J_{\lambda}}$ and hence $\sigma^{J_{\lambda}}$. {\qed}(Claim)\\

We may then finish off the Lemma as follows.
Note first that the expression in quotation marks on the right hand side,
$\eta_{\sigma} ( i )$ say, here is, if true, a member of $T^{2}_{\tau}$, being
$\Sigma^{J_{\tau}}_{2}$ in $i$. We thus shall have $\sigma \in T^{2}_{\lambda}
\lr \ex i \ul \eta_{\sigma} ( i ) \ur \in T$ and the Lemma is proven. \qed (Lemma \ref{L1.4})\\

Finally, by the last lemma the
code $x_{\lambda}$ for $J_{\lambda}$ on the scratch tape can be written by stage $\omega^{2}
\cdot \lambda + \omega + \omega$. A code for $J_{\lambda +1}$ and the diagram
$d_{\lambda +1}$ is written by stage $\omega^{2} \cdot ( \lambda +1 )$, and
$T^{2}_{\lambda +1}$ by $\omega^{2} \cdot ( \lambda +1 ) + \omega \cdot
2$. {\qed} (Lemma \ref{L1.3})\\

\section{Halting times of transfinite Turing machines}

We show that the spectrum of halting times of, {\eg}, ITTMs, is changed if we
allow multiple independent read/write heads. An extra R/W head allows for
example, $\omega_{1}^{\tmop{ck}}$, and other admissibles to become halting
times.

\subsection{Multiple Head Machines}

Multitape/multihead ITTMs are not going to be able to compute more functions
than the standard ITTMs, since the latter can simulate the former, either
by directly programming them in, or else just appealing to the observation
that multihead machines are still absolute to $L $ in their actions. So knowing
$L_{\beta}$ is to know the action of such a machine up to stage $\beta$ as a
sequence of snapshots of all its cell values, head positions, and states as a
$\Sigma_{2}$-recursion over $L_{\beta}$, and thus in turn to know its
behaviour from the theory of $L_{\beta}$.  The latter in turn the single-head machine can
discover through the Theory Machine of {\cite{FrWe08}}. The only change that
is relevant for this discussion is the change in possible halting times, which
indeed the extra heads can facilitate. Recall that no admissible ordinal is
the halting time of any $P_{e} ( n )$ for a standard single head ITTM
({\cite{HL}}). The following shows that all ordinals up to the first
$\Pi_{3}${\tmem{}}-reflecting ordinal are such halting times.

\begin{theorem}
  \label{Th2.3} (i) Suppose we consider multiple head {\tmem{ITTM}}s. (a) Let
  $\alpha$ be $\Pi_{3}${\tmem{}}-reflecting. Then $\alpha$ is not
  $\tmop{ITTM}$-clockable.
  
  (b) The least $\Pi_{3}${\tmem{}}-reflecting $\alpha_{0}$ starts the first
  gap in such clockables. In general, if $\alpha$ is not
  $\Pi^{1}_{3}${\tmem{}}-reflecting, then $\alpha$ does not start a gap; 
  hence for such $\alpha$, if it is a limit of clockables, it is itself clockable.
  
 (ii)  The same holds true for {\tmem{OTM}}s if again we allow multiple heads
  and tapes. 
\end{theorem}

{\pf} (i)(a) We assume initially that there are just three tapes as on an
ITTM, which have zeros written everywhere. Assume for a contradiction that
$P_{e} ( 0 ) \da^{\alpha}$.  In an enumeration $\pa{C^{k}_{n}}_{n}$ \ of the
cells let $C^{k}_{n} ( \eta ) =j$ indicate that the $n$'th cell of tape $k$,
at time $\eta$ has contents $j<2$ for $k<3$, where $k=0,1,2$ representing the
Scratch, Output, and Input tapes respectively. Let $R^{0} ( \tau ) \nocomma
\nocomma ,R^{1} ( \tau ) ,$be respectively, the position of two R/W heads on
the Input/Scratch and Output Tapes respectively (we might easily imagine a
third head $R^{2}$, and thus in effect we set for any $\tau$ $R^{0} ( \tau )
=R^{2} ( \tau )$); let $I ( \tau )$ be the current instruction or state number
($\nobracket <  \omega )$ at time $\tau$ for $P_{e}$, for $\tau \leq \alpha$.
We note that for $\tau < \alpha$ that $R^{k} ,I$ are
$\Sigma_{1}^{L_{\alpha}}$-definable functions of $\tau$. We shall assume that
at limit times $\lambda$, a head returns to the Liminf of its positions on its
tape if this is finite, otherwise the head returns to 0 the starting position.
Mathematically put: at time $\lambda$, $R^{k} ( \lambda ) =
\tmop{Liminf}^{\ast}_{\beta \imp \lambda} R^{k} ( \beta )$ for $k<3$, where we
define $\tmop{Liminf}^{\ast}_{\beta \rightarrow \lambda} R^{k} ( \beta ) =
\tmop{Liminf}_{\beta \rightarrow \lambda} R^{k} ( \beta )$ if the latter is less than
$\omega$, otherwise set this to $0.$

Now abbreviate $r^{k} =R^{k} ( \alpha )$ for $k=0,1$ (and set also $r^{2}=r^{0}$
whenever it appears as $r^{k}$ for $k=2$).
By our assumption, the machine halts because of a particular 
configuration under the R/W heads $C_{r^{k}}^{k} ( \alpha )$, the reason
being that it is in a particular state $I ( \alpha )$ whilst reading precisely those
cells. By simple $\Pi_{2}$ reflection at the admissible $\alpha$, $I ( \beta )
=I ( \alpha )$ for an unbounded (indeed closed) set $C$ of ordinals below
$\alpha$.

{\tmem{Case 0}} \ \ Both $r^{0} ,r^{1} >0$.

Then again by an additional $\Pi_{2}$-reflection clause there is an unbounded
set $D_{0} \sset C$ of ordinals $\gamma$ below $\alpha$ where $r^{k} =R^{k} (
\gamma )$ for both $k<2$ simultaneously. Note further that ``$C_{r^{k}}^{k} (
\alpha ) =j^{k} <2$, for $k<3$'' is expressible as a Boolean combination of
$\Sigma_{2}$ and $\Pi_{2}$ sentences. (For example, either $j^{k} =1$, and
then $C_{r^{k}}^{k} ( \tau ) =1$ if it is 1 for all sufficiently large $\tau <
\alpha$, or else it is 0, which happens by Liminf rules if we had that for all
$\tau_{0} < \alpha$, there is a $\tau > \tau_{0}$ with $C_{r^{k}}^{k} ( \tau )
=0$. So this is $\Sigma_{2} \vee \Pi_{2}$.) \ We thus have some unbounded $D
\sset D_{0}$ where these facts reflect to $\gamma \in D$. But then for such
$\gamma$ the conditions on $I ( \gamma ) ,C^{k}_{r^{k}} ( \gamma )$ were all
there to halt exactly as at stage $\alpha$. A contradiction.

{\tmem{Case 1}} \ Just one of $r^{k} =0.$

Suppose w.l.o.g. this is $r^{0}$. Then there is an unbounded $D \sset C$ of
ordinals $\gamma$ where $r^{1} =R^{1} ( \gamma )$ and $C^{1}_{r^{1}} ( \gamma
) =1$, just as in Case 0.

{\tmem{Case 1a) }} $\tmop{Liminf}_{\beta \rightarrow \alpha} R^{0} ( \beta ) =
\omega$.

Then \ \ $\all z< \omega \ex \tau_{z} \all \tau \in ( \tau_{z} , \alpha ) ( 
R^{0} ( \tau ) >z )$. \ However this is a $\Pi_{3}$-sentence and we may
reflect this down, simultaneously with the above $\Sigma_{2} \vee
\Pi_{2}$-statements about $r^{1}$, to some $\gamma \in D$. But again this
yields the conditions for halting at a stage $\gamma < \alpha$ as $R^{0} (
\gamma ) =0=r^{0}$ - a contradiction.

{\tmem{Case 1b) }} $\tmop{Liminf}_{\beta \rightarrow \alpha} R^{0} ( \beta ) =0$.

Then $\all \tau_{0} \ex \tau \in ( \tau_{0} , \alpha ) R^{0} ( \tau ) =0$.
This is $\Pi_{2}$, so again can be reflected to some $\gamma \in D$ as in Case
1a), reproducing the halting conditions at time $\gamma$. Again a
contradiction.

{\tmem{Case 2}} \ Both $r^{k} =0$.

Then the argument is just a mixture of those of the previous Case, depending
on the reason for each $r^{k}$ to be 0. We leave this as an exercise to the
reader. {\qed} Thm \ref{Th2.3} (i)(a)\\

{\nod}To continue we shall need a special case of the following Lemma.

\begin{lemma}
  \label{L2.3}Let $\beta > \omega  $ be such that $L_{\beta}$ is not
  $\Pi_{3}$-reflecting, is not a model of $\Sigma_{1}$-Separation, but is an
  admissible model of ``$V= \tmop{HC}$''. Then there is a $( \Pi_{2} \wedge
  \Sigma^{}_{2} )^{L_{\beta}} ( \{ p \} )$ function $G: \omega \imp \beta$
  which is cofinal, where $p=p^{1}_{\beta}$ is the standard parameter of
  $L_{\beta}$, but so that for any $\beta' < \beta  $ with $p=p^{1}_{\beta'
  ,}$ $G^{L_{\beta'}}$ is not a total function.
\end{lemma}

{\pf}Let $\prec -,- \succ : \omega \times \omega \imp \omega$ be a recursive
pairing function. Our assumptions imply that $h^{1}_{\beta}$``$\omega \times (
\omega \cup \{ p \} ) =L_{\beta}$. (Here $h=h^{1}_{\beta}$ is the standard -
and uniformly defined - $\Sigma_{1}$-skolem function for all limit
$L_{\beta}$). Suppose $\varphi [ x ] \equiv \all u \ex v \all w \psi ( u,v,w,x
)$ for some $x \in L_{\beta}$ is a non-reflecting $\Pi_{3}$-formula, chosen so
that $L_{\beta} \models \varphi [ x ]$. This is equivalent to:

$$(+)\quad  \all \prec n,i \succ \,\,\,\,\, < \,\,\, \omega   \ex \alpha \left[ h ( i,n,p ) \ua \vee
\ex v \in L_{\alpha}   \right. \left. \all w \psi ( h ( i,n,p ) ,v,w,x )
\right]$$

{\nod}holding in $L_{\beta}$. Let 
$$
\begin{array}{lcl}
G ( \prec n,i \succ ) & = & \mbox{ the least } \alpha \mbox{
such that }\ex v \in L_{\alpha} \all w \psi ( h ( i,n,p ) ,v,w,x ), \mbox{ if }h ( i,n,p )
\da \\
& = & \,\, 0,  \mbox{ if }h ( i,n,p ) \ua .
\end{array}
$$
\nod Then $G$ is $( \Pi_{2} \wedge \Sigma^{}_{2} )^{L_{\beta}} \{ p \}$ - the
`leastness of $\alpha$' requirement contributing a $\Pi^{L_{\beta}}_{2}$
clause:
$$\all \alpha' < \alpha
\all v \in L_{\alpha'}    \neg \all w \psi ( \nobracket \nobracket ( h (
i,n,p ) ,v,w,x ) ).$$
\nod We claim that for any limit $\beta' < \beta$ with $x \in L_{\beta'}$, and
$p=p^{1}_{\beta'}$ (and there are unboundedly many such $\beta'$ below
$\beta$), that $G^{L_{\beta'}}$ is not a total function. Our assumptions show
that again \ $h^{1}_{\beta'}$``$\omega \times ( \omega \cup \{ p \} )
=L_{\beta'}$ (although for those $h^{1}_{\beta} ( i,n,p ) \in L_{\beta} \back
L_{\beta'}$ we shall have $h^{1}_{\beta'} ( i,n,p ) \ua$). So were \
$G^{L_{\beta'}}$ total, we should have that $( \varphi [ x ] )^{L_{\beta'}}$
which is impossible.{\qed}Lemma \ref{L2.3}.\\

Now for limit $\beta' \leq \beta  $ let $n ( \beta' ) \dfs \tmop{maximal}  n'
\leq \omega$ s.t. $L_{\beta'} \models$``$  \all \prec \bar{n} ,i \succ \leq n'
\left( \prec \bar{n} ,i \succ \in \dom ( G ) \right)$''. Then for
appropriate $\beta' < \beta$ we must have $n ( \beta' ) < \omega$, whilst $n (
\beta ) = \omega$.

\begin{lemma}\label{L2.3.5} $\tmop{Liminf}_{\beta' \rightarrow \beta}  n ( \beta' ) =
\omega$, and 
$G ( \prec
\bar{n} ,i \succ )^{L_{\beta'}}\imp G (
\prec \bar{n} ,i \succ )^{L_{\beta}}$.
\end{lemma}

{\pf} As $\beta' \rightarrow \beta$ the value of $G ( \nobracket \nobracket \prec
\bar{n} ,i \succ )^{L_{\beta'}}$ eventually settles to the final value $G (
\nobracket \nobracket \prec \bar{n} ,i \succ )^{L_{\beta}}$. To see this note that for a
given $\beta' < \beta$ the fact that $\all v \in L_{\beta'}   \neg \all w \psi
( h ( i,n,p ) ,v,w,x )^{L_{\beta}}$ is a $\Sigma_{1} \{ p,x \}^{L_{\beta}}$
fact by admissibility of $\beta$, and so if true in $L_{\beta}$ will also be true in sufficiently
large $\beta' < \beta$. Hence larger and larger initial portions of the domain
of $G$ become correctly calculated as $\beta'$ increases. \\ \mb \hfill {\qed} Lemma \ref{L2.3.5}

\begin{lemma}
  \label{L2.4}Suppose that $\omega < \beta$ is not $\Pi_{3}$-reflecting and is
  a supremum of halting times of \ 2-head {\tmem{ITTM}} computations on
  integer input. Then $\beta$ itself is such a halting time; it thus does not
  start a gap in such {\tmem{ITTM}} clockables.
\end{lemma}

{\pf}We in essence run the ``Theory Machine'' \ that computes successively the
constructible levels of the $L_{\tau}$ hierarchy, for $\tau < \beta < \Sigma$;
by time $\beta$ this machine has written (and overwritten) codes successively
for all the $L_{\tau}$. We assume that $\beta$ is admissible, since if not
then we already know that $\beta$ is clockable, as only admissibles start gaps
even for single head machines (see {\cite{W09}}, Thm. 50). Let $x \in
L_{\beta}$ and $\varphi [ x ] \equiv \all u \ex v \all w \psi ( u,v,w,x )$. We
assume $\varphi [ x ]$ fails to be true in any $L_{\tau}$ for $\tau < \beta$
but is true in $L_{\beta}$. \ Let $G: \omega \imp L_{\beta}$ be a cofinal,
$\Pi_{2} \wedge \Sigma^{L_{\beta}}_{2}$ function defined as above (parameter
free), which arises from the fact that $L_{\beta} \models$``$V= \tmop{HC}$''
and is not a model of $\Sigma_{1}$-Sep. (This latter requirement and the
parameter free-ness of the definition follows from $\beta$ being a supremum of
halting times.) We may
thus assume without loss of generality that $x \sset \omega$.

Now consider the amendment to the Th{\nod}eory Machine that is also recording
on a separate tape $T_{1}$ with its own R/W head (this can be the Output Tape
of the Theory Machine,
$\la C_{n}^{1} \ra_{n}$, and we can arrange that this tape is otherwise not
used by the Theory Machine) an initial sequence of: the first cell $C^{1}_{0}$ which
is kept at 0, then in the next $n ( \beta' ) +1 $ cells 1's are written with
the head $R^{1}$ resting on the largest cell with a 1 after this write, namely
$C^{1}_{n ( \beta' ) +1}$, and thereafter by $0$'s, to indicate, that at the
current stage $\beta' \leq \beta$, $n ( \beta' )$ is our maximal $n \leq   \,
\omega  $ as defined above. The value of $n ( \beta' )  $ of course being
found from inspection of the $\Sigma_{2}$-theory of $L_{\beta'}$ that the TM
provides.

Then $\beta' \rightarrowtail n ( \beta' )$ is a, not necessarily monotone,
function into $\omega$, but with $\tmop{Liminf}_{\beta' \rightarrow \beta}  n
( \beta' ) = \omega$ (see Lemma \ref{L2.3.5}). Thus at time $\beta$
the R/W head of the tape $T_{1}$, having been sent out to $\omega$, is returned to its $0$'th cell $C^{1}_{0} :R^{1} (
\beta ) =0$. Thus the head is reading a 0, and this is the first time this
happens at any stage $\beta' >0$. So we can program a halt at this
stage.{\qed} Lemma \ref{L2.4} \& Thm. \ref{Th2.3}(i)(b).\\

{\pf}of Theorem \ref{Th2.3} (ii)(a).

We adapt the above arguments for when we have an OTM with multiple heads, let
us take three independent heads on three  tapes, with the same notations of
$R^{k}_{} ( \tau ) \in \tmop{On} ,C^{k}_{\delta} ( \tau ) <2$ (now for the
$\delta$'th cell of tape $T^{k} \dfs \pa{C_{\delta}^{k}}_{\delta \in
\tmop{On}}$) for $k<3$ and $I ( \tau )$ as before.

(a) We assume that $\alpha$ is $\Pi_{3}$-reflecting, and show that it is not
$\tmop{OTM}$-clockable under this regime. So assume for a contradiction that
$P_{e} ( 0 ) \da^{\alpha}$. Now the R/W heads can be at various
combinations of positions once more. Let $r^{k} =R^{k} ( \alpha )$ for $k<3$.
We consider the case that $r^{1} \in \tmop{Lim}$ and that \ $r^{0} ,r^{2}$ are
successor ordinals. This is the simplest but paradigm case: with $r^{0}
,r^{2}$ successor ordinals, there is by $\Pi_{2}$-reflection as above a cub $C
\sset \alpha$ where the heads on these two tapes $T^{0} ,T^{2}$ are at the
same position (thus $R^{k} ( \tau ) =r^{k}$ for $k=0,2$), and the instruction
number $I ( \tau )  $ is the same, for any $\tau \in C.$ We then just
concentrate on $r^{1} \in \tmop{Lim}$. Let us further abbreviate $r^{1}$ as
$\rho$.

As $\rho >0$ there is some $\tau_{0} < \alpha$ so that for no later stage
$\tau \geq \tau_{0}$ can we have $R^{1} ( \tau ) \in \tmop{Lim}$ and $R^{1} ( \tau
+1 ) <R^{1} ( \tau )$ (``no further going left at limit cells'').\\

{\nod}{\tmem{Case 1}} $\ex \tau_{0} < \alpha \all \tau \in ( \tau_{0} , \alpha
)  R^{1} ( \tau ) < \rho$.\\

If $\rho < \alpha$ then we can assume that $\rho < \tau_{0}$. In any case, by
definition $\rho = \tmop{Liminf}_{\tau \rightarrow \alpha} R^{1} ( \tau )$, so
we have that $\rho$ under our case hypothesis can be defined as:
$$y= \rho
\equi \all z<y \ex \tau_{z} \all \tau \in ( \tau_{z} , \alpha ) ( y>  R^{1} ( \tau
) >z).$$

{\nod}By $\Pi_{3}$-reflection there is an admissible $\beta < \alpha$ with
$\tau_{0} < \beta \wedge \all \tau \in ( \tau_{0} , \beta )  R^{1} ( \tau ) <
\rho$ and also
$$ \all z< \rho \ex \tau_{z} \all \tau \in ( \tau_{z} , \beta ) (\rho >  R^{1}
( \tau ) >z ).$$

{\nod}(In the case that $\rho = \alpha$, then the $\rho$ in both these last
expressions is then taken to be $\beta$.)\\

{\nod}{\tmem{Case 2}} Otherwise (which implies $\rho < \alpha  $ as $R^{1} (
\tau ) < \alpha$ for any $\tau < \alpha$).\\

By increasing $\tau_{0}$ if need be in this case we may assume $\all \tau_{2} \in (
\tau_{0} , \alpha )  R^{1} ( \tau_{2} ) \geq \rho$ (as there is no going left
at limits beyond $\tau_{0}$).

We now argue that the conditions for $P_{e} ( m ) \da$ at time $\alpha$ also
occurred at an earlier time $\bar{\alpha}$. We've remarked above that we can
assume by $\Pi_{2}$ reflection at the admissible $\alpha$ that there is a
closed unbounded set $C \sset \alpha$ of $\bar{\alpha}$, \ with $I (
\bar{\alpha} ) =I ( \alpha ) =s$, and $R^{k} ( \bar{\alpha} ) =r^{k}$ for
$k=0,2$.

In {\tmem{Case 1}} we have that $R^{1} ( \beta ) = \rho$ (respectively $R^{1}
( \beta ) = \beta$ if $\nobracket \rho = \alpha )$. By the case hypothesis the
value of the cell $C_{\rho} ( \tau )$ is unchanged in the interval $( \tau_{0}
, \alpha )$ and so also in $( \tau_{0} , \beta )$. \ We may also assume by an
additional reflection clause, that $\beta$ is in the $\Pi_{2}$-definable set
$C \nosymbol .$ We thus have that cell contents of $C^{k}_{\rho}$, head
positions $R^{k} ( \beta )$, and state number $I ( \beta ) =s$ are the same at
times $\beta$ and $\alpha$. Thus all the ingredients are there for $P_{e}$ to
halt at time $\beta < \alpha$.

In {\tmem{Case 2}} knowing that $\all \tau_{2} \in ( \tau_{0} , \alpha ) 
R^{1} ( \tau_{2} ) \geq \rho$, then the $\Pi_{2}$ fact about $\alpha$ (i) that
$R^{1} ( \alpha ) = \rho$, and additionally (ii) $C^{1}_{\rho} ( \alpha ) =0$
(if the latter holds) reflect to an unbounded (in fact closed) subset $E$ of
$C$ above $\rho$. If $C^{1}_{\rho} ( \alpha ) =1$ then we have $C^{1}_{\rho} (
\tau' ) =1$ for a tail of $\tau' < \alpha$, and so also on a tail of $\tau'
\in E_{0}$ where $E_{0}$ reflects (i) alone. Again the conditions for halting at
$\alpha$ are reflected to any $\bar{\alpha} \in E$ (or respectively $E_{0}$),
and so $P_{e}$ must halt earlier. {\qed} Theorem \ref{Th2.3} (ii)(a)\\

The following discussion yields the version of Theorem \ref{Th2.3} (ii)(b) for
OTMs.

\begin{definition}
  $\sigma_{1}$ is the first {\tmem{stable}} ordinal: that is, it is the least
  $\sigma$ so that $L_{\sigma} \prec_{\Sigma_{1}} V$.
\end{definition}

\begin{remark}
  Now if any \tmem{OTM} program halts on some integer input, then this $\Sigma_{1}$
  fact must reflect to $L_{\sigma_{1}}$, hence this halt must be at a stage
  $\alpha < \sigma_{1}$, and moreover with tape contents an element of
  $L_{\sigma_{1}}$. We note that halting times of \tmem{OTM}s on integer input, are
  cofinal in $\sigma_{1}$: by the definition of $\sigma_{1}$, there are
  arbitrarily large $\alpha < \sigma_{1}$ where some new
  $\Sigma_{1}$-sentence, $\varphi$ say, is first true at $L_{\alpha +1}$, but
  not at $L_{\alpha}$. As Koepke shows we may run an \tmem{OTM} that produces codes
  for $L_{\beta}$'s. Then we may adopt this to halt when it finds one in which $\varphi$
  is true; this will take more than $\alpha$, but less than $\sigma_{1}$, many
  steps.
\end{remark}

\begin{lemma}
  \label{L2.8} Suppose that $\omega < \beta$ is not $\Pi_{3}$-reflecting and
  is a supremum of halting times of multi-head-\tmem{OTM} computations on integer
  input. Then $\beta$ itself is such a halting time; it thus does not start a
  gap in the multihead-\tmem{OTM} clockables.
\end{lemma}

{\pf} This proof should just be: run the same construction as Lemma \ref{L2.4}.
We need to make some observations and adjustments pertinent to the
OTM-scenario. Firstly, because $\beta$ is assumed to be a supremum of such
OTM-halting times, we have that $L_{\beta} \models$``$V= \tmop{HC}$''. This is
because, if $P_{e} ( k ) \da^{\alpha}$ then this is a $\Sigma_{1}$-fact, first
true in $L_{\alpha +1}$, and it is then easily seen that $\alpha$ is not a
cardinal in $L_{\alpha + \omega}$. (Otherwise, by the Levy reflection Theorem
applied in $L_{\beta}$, every such $\Sigma_{1}$-fact would be true at a level
before $\sigma^{L_{\beta}} < \omega_{1}^{L_{\beta}} < \beta$ by assumption.)
Hence any such $L_{\alpha + \omega} \models$``$\omega$ is the largest
cardinal''. (And $\alpha < \sigma_{1}$.) As this holds for unboundedly many
$\alpha + \omega$ below $\beta$, it holds for $L_{\beta}$.

Secondly, we assume that $\beta$ is admissible. (If $\beta$ is inadmissible,
then known methods can find a suitable program $P_{e} \da^{\beta}$. We leave
this to the reader.)

Essentially we run the same proof: here however we run Koepke's version of
the Theory Machine: we assume on the scratch tape that (adapting the truth
predicate machine of {\cite{K05}}) that the program computes successively
the constructible levels of the $L_{\tau}$ hierarchy, now for $\tau < \beta$;
by the admissible time $\beta$ this machine has written codes successively for all the $\tau
< \beta$. As before we let $\varphi [ x ] \equiv \all u \ex v \all w \psi (
u,v,w,x )$, where we assume $\varphi [ x ]$ fails to be true in any $L_{\tau}$
for $\tau < \beta$, but is true at $L_{\beta}$. \ Define $G$ as before: then still $G: \omega \imp \beta$
is a cofinal function, $( \Pi_{2} \wedge \Sigma_{2} )^{L_{\beta}}$ (parameter
free). That $G$ exists uses the first observation that \ $L_{\beta}
\models$``$V= \tmop{HC}$'' and again, being a supremum of halting times, there
are no proper $\Sigma_{1}$-substructures of $L_{\beta}$ and hence
$p^{1}_{\beta} = \emp$. Now run the argument exactly as before. The only
difference is that our definition of OTM behaviour has that the R/W head on
tape $T^{k }$ does not return to $C^{k}_{0}$ at some limit stage by
overshooting the end of the tape, but can do so at a successor stage, either
in the natural course of the Turing program or because it attempted to
``move left whilst on a limit cell''. But all is well and at the very end of
the above proof, at time $\beta$ the R/W head of $T^{1}$ is not returned to
its $0$'th cell $C^{1}_{0}$, but is placed on $C^{1}_{\omega}$ where however
it is reading a 0, and again this is the first time this happens at any stage
$\beta' >0$. So again we can program a halt at this stage. {\qed} Lemma \ref{L2.8} \& Thm. \ref{Th2.3} (ii).

\begin{corollary}
  \label{C2.8}If $\gamma$ is a supremum of halting times in either the \tmem{ITTM} or
  \tmem{OTM} multihead models, then $\gamma$ starts a gap iff $\gamma$ is
  $\Pi_{3}$-reflecting.
\end{corollary}

\subsection{Single Head Machines}

For single head OTMs the situation is more subtle: there are some admissible
ordinals that are not $\Pi_{3}$-reflecting but which nevertheless start gaps.
One may even ask:\\

\nod{\bf{Question}}: Is $\omega_{1}^{\tmop{ck}}$ clockable by a one head OTM
program?\\

Note that the ITTM argument of {\cite{HL}} which shows that
$\omega_{1}^{\tmop{ck}}$ is not clockable does not apply, since although there
the machine has one head, we have diverged from their limit rules. The
argument above, directly appealing to the non-$\Pi_{3}$-reflection at
$\omega_{1}^{\tmop{ck}}$, or alternatively that of Carl {\cite{Carl20}} where
ITRMs (which do clock the $\omega_{n}^{\tmop{ck}}$'s) are simulated on OTMs,
both are with multiple heads and tapes, and thus do not apply either. We shall
answer this question below.

 As multihead machines subsume single head ones we
have, as a corollary to Thm.\ref{Th2.3}(ii):

\begin{corollary}
  Let $\gamma$ be $\Pi_{3}$-reflecting. Then $\gamma$ is not clockable by a
  single head $\tmop{OTM}$. Hence such $\gamma$ either starts a gap, or are
  interior to one.
\end{corollary}

We now characterise the rest of the ordinals that start gaps. Some admissible ordinals that are non-$\Pi_{3}$-reflecting are clockable, some are not. 
Let $E^{\ast}$ be the set of limit points of $E$ the class of
$\Sigma_{2}$-extendibles.

\begin{theorem}
  Let $\gamma$ be admissible, but not $\Pi_{3}$-reflecting. Suppose $\gamma$
  is a supremum of halting times for single head {\tmem{OTM}}s. Then
  
  \ \ \ \ \ \ \ \ \ \ \ \ \ \ \ \ \ \ \ \ \ \ \ \ \ \ \ \ \ \ \ \ $\gamma$ is
  clockable $\equi   \gamma$ is not a limit of $\Sigma_{2}$-extendibles.
  
  {\nod}Hence the class of such $\gamma$ which start gaps consists of
  precisely those $\gamma$ in $E^{\ast}$.
\end{theorem}

{\pf} Let $\gamma$ be as assumed: admissible, but not $\Pi_{3}$-reflecting,
and a supremum of halting times for single head {OTM}s.
Consider, as for ITTMs, the example case of an {OTM} with a single R/W head at
position $R ( \tau )$ at time $\tau$, but three tapes enumerated
$\pa{C^{k}_{\alpha}}_{\alpha \in \tmop{On}}$ for $k<3$ with content at time $\tau$
$\pa{C^{k}_{\alpha} ( \tau )}_{\alpha \in \tmop{On}}$. As for ITTMs we
envisage the head reading at any moment simultaneously a triple
$\pa{C^{k}_{\alpha}}_{k<3}$.

$\left( \imp \right)$ We assume that $\gamma$ is both a limit of extendibles
and is clockable for a contradiction. At time $\gamma$ let $I ( \gamma ) =i<
\omega$. As above by $\Pi_{2}$-reflection there is a closed and unbounded $C
\sset \gamma$ with $I ( \delta ) =i$ for $\delta \in C$. In the argument that follows, the main case is when the position of the head at time $\gamma$, here called $\rho$ once more, is at neither extreme of $0$ or $\gamma$. The dynamics of head movement allied with the limit rules will ensure this $\rho$ can not be a successor ordinal; coupled wth the admissibility of $\gamma$, it cannot be a limit of limit ordinals. We then shall deduce that from some point onwards it is behaving isomorphically to an ITTM, in that it is restricted to an $\omega$-sequence of cells $ [ \bar{\rho} ,\bar{\rho} + \omega )$ where $\bar{\rho} + \omega =\rho $. Then the assumption that $\gamma$ is a limit of extendibles plays its role.

{\nod}{\tmem{Case 1}} $\rho = \gamma$.

Then the R/W head is reading the three cells $C^{k}_{\gamma} ( \gamma ) =0$
for each of $k<3$ for the first time. This is $\Pi_{2}$-expressible (``$\all
\delta \ex \tau R ( \tau ) \geq \delta$'') and so reflects to an $\alpha \in
C$, with $C^{k}_{\alpha} ( \alpha ) =0$ ($\nobracket k<3 )$ and the conditions
for halting are there at time $\alpha  $ - a contradiction.

{\nod}{\tmem{Case 2}} $0< \rho < \gamma$.

Then there is some $\tau_{0} < \gamma$ so that for all $\tau \in ( \tau_{0} ,
\gamma )$, if $R ( \tau ) \in \tmop{Lim}$, then $R ( \tau +1 ) \not{<} R (
\tau )$ (``No going left at a limit''). In particular this implies that we
have just the two alternatives:

{\em Either:} $\ex \tau_{0} < \gamma \all \tau \in ( \tau_{0} , \gamma ) R ( \tau )
\in [ \rho , \rho + \omega )$. By an additional $\Pi_{2}$-reflection clause
then (using $\tau_{0}$ as a parameter for example), there is a closed and
unbounded set $D \sset C \sset \gamma$ with $\tau \in D \imp R ( \tau ) =
\rho$. Similarly there is an unbounded set $B \sset D$ to which the pattern \
$C^{k}_{\rho} ( \gamma )$ for $k<3$ reflects as $C^{k}_{\rho} ( \alpha )$ for $\alpha\in B$, and for which the conditions for halting
are present. Contradiction.

{\em Or:} $\ex \tau_{0} < \gamma \all \tau \in ( \tau_{0} , \gamma ) R ( \tau ) <
\rho$. \ Note $\rho$ cannot be a successor ordinal. Nor can it be a compound
limit. For, let $\rho_{n} < \rho_{n+1}$ be an increasing sequence of limit
ordinals with supremum $\rho$ and with $\pa{\rho_{n}}_{n} \in L_{\gamma}$. As
$R ( \gamma ) = \rho$, for each $n$ there is a least $\tau_{n} > \tau_{0} ,
\rho$ so that $R ( \tau_{n} ) \geq \rho_{n}$. On the one hand $n
\rightarrowtail \tau_{n}$ is a total $\Sigma_{1}^{L_{\gamma}} \left( \left\{
\pa{\rho_{n}}_{n} , \tau_{0} \right\} \right)$ definable function, but on the
other if $\tilde{\tau} =  \sup_{n}   \tau_{n}$, then $R ( \tilde{\tau} ) =
\rho$ (here we are implicitly using that $\rho >0$ to say $\tau' > \tau_{n}
\imp R ( \tau' ) \geq R ( \tau_{n} )$). However our case hypothesis implies
that $\tilde{\tau} = \gamma$ - a contradiction to the admissibility of $\gamma$.

We are thus left with $\rho = \bar{\rho} + \omega$ with $\tmop{Lim} (
\bar{\rho} )$. Then $\ex \tau_{0} < \gamma \all \tau \in ( \tau_{0} , \gamma )
R ( \tau ) \in [ \bar{\rho} , \bar{\rho} + \omega )$. However then the
program and OTM are behaving, within this interval of time, exactly as for
an ITTM (with Liminf rules for the position of its R/W head of course) but
working on the $\omega$ sequence of cells $[ \bar{\rho} , \bar{\rho} + \omega
)$ rather than $[ 0, \omega )$. But now if $\xi > \tau_{0}, \rho$ is the least
$\Sigma_{2}$-extendible, with $L_{\xi} \prec_{\Sigma_{2}} L_{\sigma}$, this
version of an ITTM will enter a final loop at time $\xi$ (as halting at a time
$\gamma' < \xi < \gamma$ is ruled out). But $\gamma$ is a limit of
$\Sigma_{2}$-extendibles, so $\xi < \sigma < \gamma$. (If $\rho < \xi < \gamma
< \sigma$, as $\gamma$ is posited a limit of extendibles there would be a
further  $\rho < \xi < \xi' < \sigma' < \sigma$ and $\xi' < \gamma$. But by
$\Sigma_{1}$-reflection applied in $L_{\sigma}$ to the pair $( \xi' , \sigma'
)$ there would be in any case another extendible pair $( \xi'' , \sigma'' )$
in $L_{\xi}$ with $\rho < \xi'' < \sigma'' < \gamma$ as required.) \ So in
fact the OTM will also loop rather than halt. Another contradiction.

{\nod}{\tmem{Case 3}} $\rho =0$.

If $\ex \tau_{0} < \gamma \all \tau \in ( \tau_{0} , \gamma ) R ( \tau )
\in [ 0, \omega )$ then this is argued as in Case 2 (indeed it is just another
example of the `Either' clause of Case 2). But $\rho =R ( \gamma )$ could be
$0$ because for arbitrarily large $\delta < \gamma$ we ``went left at a
limit'': $R ( \delta ) \in \tmop{Lim}$, \ $0=R ( \delta +1 ) <R ( \delta )$.
This is $\Pi_{2}$-expressible. By the Liminf rule this means there is a{\tmem{
closed}} and unbounded set $C$ of $\delta$ below $\gamma$ with $R ( \delta ) =0$.
So again by reflection there is an earlier  $\tau \in C$ with $R ( \tau ) =0$, and
arguing as above, with  pattern match-up: $C^{k}_{0} ( \gamma ) =C^{k}_{0} ( \tau )$ for $k<3$.
Again we have the conditions for halting at this $\tau$, a contradiction.
{\qed}$\left( \imp \right)$\\

$( \longleftarrow )$ Assume $\gamma$ is not a limit of extendibles and we
shall show $\gamma$ is clockable. Let $\delta_{0}$ be the largest extendible,
or limit of such, less than $\gamma$. As $\gamma$ is a limit of clockables,
there is $\delta \in ( \delta_{0} , \gamma )$ and some $f,k$ so that $P_{f} (
k ) \da^{\delta}$. This is a new $\Sigma_{1}$-fact true in $L_{\delta +1}$. \
This implies by standard constructibility theory, that $L_{\delta +2}
\models$``$V= \tmop{HC}$''. Let $x_{\delta} \in L_{\delta +2} \cap \tmop{WO}$
be the $L$-least code for $\delta$. We describe a single head 3-tape OTM
procedure.

(1) The Koepke $L$-machine runs constructing levels of $L$ until it is seen
that there is a level of $L$ witnessing $P_{f} ( k ) \da^{\delta}$, and then
$x_{\delta}$ is found. The discovery of $x_{\delta}$ is a signal to advance to
the rest of the main process. \ From this point onwards, thus some $\tau_{0}$
onwards, $R ( \tau ) \in [ 0, \omega )$, {\ie}the process will be confined to
just the first $\omega$ many cells of each tape. Only at time $\gamma$ shall we
 have $R ( \gamma ) = \omega = \tmop{Liminf}_{\tau \rightarrow \gamma} R
( \tau )$. At this point there will be flags, in other words a cell
configuration, in $\pa{C^{k}_{0}}_{k<3}$, which we can place there now, which
will program an immediate halt, thus demonstrating the clockability of
$\gamma$. Also at this point $\tau_{0}$ we can write $x_{\delta}$ to the  first
$\omega$ cells of the input tape, whilst those of the Scratch and Output tapes
are reset to zeroes. After these preliminaries we proceed to:

(2) apply the notation of Lemma \ref{L2.3} above, which provided a cofinal
function $G: \omega \rightarrow \gamma$ which was $( \Pi_{2} \wedge
\Sigma^{}_{2} )^{L_{\gamma}}$ but so that for any $\beta < \gamma $,
$G^{L_{\beta}}$ was not a total function. (We note that no parameters are
needed here, as by our assumption $\gamma$ is the supremum of clockables
$\gamma'$, and
so the supremum of ordinals where new $\Sigma_{1}$-facts become true. If $p^{1}_{\gamma}\neq
\emp$ then, by standard arguments, we have that $L_{\max p^{1}_{\gamma}}\prec_{\Sigma_{1}}L_{\gamma}$ and so there are no new $\Sigma_{1}$-facts becoming true in the interval $(\max p^{1}_{\gamma}, \gamma)$.) 
We let $n ( \alpha )$
be as defined after that Lemma. We ensure whilst $n ( \alpha ) =n_{0}$ say,
that the program does not visit any cell $C^{k}_{t}$ for any $t<n_{0}$.
Since $\tmop{Liminf}_{\alpha \rightarrow \gamma} n ( \alpha ) = \omega$, this
will ensure $R ( \gamma ) = \omega$ as desired.

(3) In order to satisfy the workspace imposed by the insistence that we only use cells in $ [ n_{0} ,
\omega )$, we divide up the scratch tape of our master program {\via}a (1-1)
(ordinary) recursive function $\pi : \omega \times 4 \times \omega \imp
\tmop{Evens}$ with the property that (i) \ $n< \min \{ \pi ( k,l,n ) \mid  k<
\omega ,l<4 \}$ and (ii) $k<k' \imp \pi ( k,l,n ) < \pi ( k' ,l,n )$. The
reason for this choice is that we shall simulate an $\omega$-sequence of
`virtual machines', $M^{n}$, each with 4 tapes, thus: $\left(
\pa{C_{k}^{l}}_{l<4,k< \omega} \right)^{M^{n}}$. Here `cell' $( C_{k}^{l}
)^{M^{n}}$ of the simulation occupies cell $C_{\pi ( k,l,n )}^{1}$ (in
particular by (i) the `tapes' of $M^{n}$ thus occupy real cells $C^{1}_{p}$
where $p>n$). That $\ran  \, \pi \sset \tmop{Evens}$ is to leave the Odds
$\pa{C^{1}_{2n+1}}_{n<\omega}$ free for scratch tape work of the master process.
(Four tapes here is clearly unnecessary: we could recursively amalgamate the
fourth tape and the output tape of $M^{n}$ and thus stick just with the usual
three tapes.)

(4) Typically a process may involve searching through, or copying of, all the
cells of an $M^{n}$ virtual tape. This could result in a Liminf value of the
actual R/W head being an unwanted $\omega$ before its desired time. To obviate
this we make otherwise redundant `insurance moves': after a fixed number, say
$10^{6}$, of steps in this searching process being run by $M^{n}$, the master
process places a mark and current instruction number of the simulated process
of $M^{n}$ at its current position, and returns the R/W head to the first
quadruple of cells $\left( \pa{C_{0}^{k}}_{k<4} \right)^{M^{n}}$ that is to
the initial group of virtual cells at the beginning of $M^{n}$'s virtual tape.
As noted at the end of (3), this still keeps the actual R/W head above the
actual cells $\pa{C^{k}_{n}}_{k<3}$. After doing this run down, it runs back
up to the mark and retrieves the next instruction number and does another
$10^{6}$ steps of the search process on the virtual machine $M^{n}$ {\etc}
This way after $\omega$ steps of its task the R/W head is not at the
$\omega$'th position. We shall refer to this as `making insurance moves'
without further specification (or indeed mention).\\

\nod{\tmem{The overall process}}\\

Idea: We define essentially a program so that with input the code
$x_{\delta}$ for $\delta$, and using Liminf rules for Instruction numbers and
R/W head positions (so once more to emphasise, {\tmem{not}} those rules of {\cite{HL}}), it will
satisfy $\tmop{Liminf}_{\beta \rightarrow \gamma'} R ( \beta ) = \omega$ first
when $\gamma' = \gamma$. (In essence one could regard this as an ITTM
program.)

(I) Let us name one of the 4 of $M^{n}$'s `tapes' $S^{n}_{}$ - the other three
being as usual for an OTM. On the input tapes of each of the $M^{n}$ we write
a copy of $x_{\delta}$. We shall run the ITTM theory machine TM from time to
time on differing $M^{n}$ using as input $x_{\delta}$ as a base. This builds
for us the $L_{\tau} [ x_{\delta} ]$-hierarchy over the real $x_{\delta}$ in
each $M^{n}$. For mildly closed limit $\tau > \delta$ $L_{\tau} [ x_{\delta} ]
\models$``$V=L$'' and indeed $L_{\tau} [ x_{\delta} ] =L_{\tau}$. We shall
reserve $S^{n}$ for accumulating G\"odel codes of $\Sigma_{1}$-sentences true
in such $L_{\tau}$. This theory we denote $T^{1}_{\tau}$.

(II) We set $n ( 0 ) =0$ and on $M^{0}$ we initiate $\tmop{TM}$ as indicated
on the content of its `input tape'. We run this up to some $\tau$ where
$\tmop{Lim} ( \tau )$ and$ $ with $L_{\tau} [ x_{\delta} ] =L_{\tau}$ and
which has a theory  $T^{1}_{\tau}$ with new $\Sigma_{1}$-sentences appearing unboundedly in $\tau$; whilst this
happens we also write at the same time $T^{1}_{\tau}$ to $S^{0}_{1}$. (All
this involves appropriate insurance moves.) The Turing jump of
$T^{1}_{\tau}$,$( T_{\tau}^{1} )'$ is recursively isomorphic to the complete
$\Sigma_{2}$-theory of $L_{\tau}$, $T^{2}_{\tau}$. $M^{0}$ computes this, and
ascertains $n ( \tau )$ from it. The master program then copies $T^{1}_{\tau}$ from
$S^{0}$ \ to $S^{n ( \tau )}$ - again using insurance moves. This describes
`Round 0'.

(III) \ Assume that at the beginning of Round $\eta$ ($\eta < \gamma$) the R/W
head is positioned at the beginning of the $M^{n ( \eta )}$ `tape', {\ie}on
the first of the quadruple of cells, so on $( C^{0}_{0} )^{M^{n ( \eta )}}$ of
$\left( \la C^{k}_{0} \ra_{k<4} \right)^{M^{n ( \eta )}}$.

(A) \ $M^{n ( \eta )}$ then runs a copy of the TM from its input tape,
$x_{\delta}$, until it sees all the $\Sigma_{1}$-sentences on $S^{n ( \eta )}$
instantiated; it continues to run $\tmop{TM}$ until it reaches a further limit
$\tau$ satisfying that new sentences in $T^{1}_{\tau}$ appears unboundedly in
$\tau$. These new sentences are added to $S^{n ( \eta )}$ as they appear. Thus
$S^{n ( \eta )}$ becomes $T^{1}_{\tau}$. From this we compute a possibly new value $n (
\eta +1 )$ (as at (II)).

(B) The master process copies $S^{n ( \eta )}$ to $S^{n ( \eta +1 )}$. This
requires insurance moves; but after this is done the R/W head is placed on
$( C^{0}_{0} )^{M^{n ( \eta +1 )}}$ of 
$\left( \left. \la C^{k}_{0} \right. \ra_{k<4} \right)^{M^{n ( \eta +1 )}}$
and the procedure returns to (A).

{\tmem{Note:}} At no point during (A) does the R/W head drop lower down on the
master tape below the cell which is $( C^{0}_{0} )^{M^{n ( \eta )}}$, and
during (B) below the cell which is $( C^{0}_{0} )^{M^{n ( \eta +1 )}}$ if $n (
\eta +1 ) <n ( \eta )$.

(IV) After a limit $\lambda  <  \gamma$ number of Rounds through (A) and (B),
as $\liminf_{\eta \rightarrow \lambda} n ( \eta )$ is still some finite number
$n_{1}$ say, we shall have (a) a theory $T^{1}_{\tau}$ written on
$S^{n_{1}}_{1}  $; and (b) by Liminf and insurance moves, the R/W head on $(
C^{0}_{0} )^{M^{n_{1}}}$; so we set $n ( \lambda ) =n_{1}$ and return to (A).

As $\eta \rightarrow \gamma$ more of the values $G^{L_{\eta}} ( n )$ reach
their final value of $G ( n ) \dfs G^{L_{\gamma}} ( n )$, and so by Lemma
\ref{L2.3} and the comment following, $\liminf_{\eta \rightarrow \gamma} n (
\eta ) = \omega$. This delivers $R ( \gamma ) = \omega$ as desired.

{\qed} ($( \longleftarrow )$ and Theorem)\\

Thus $\omega_{1}^{\tmop{ck}}$ is, in particular, the halting time of a single head OTM. We
leave the reader to verify that an inadmissible supremum of halting times is
also a halting time. More generally we then have:

\begin{corollary}
  \label{C2.11} For single head {\tmem{OTM}}s let $\gamma$ be a supremum of
  halting times. Then:
  
  \ \ \ \ \ \ \ $\gamma$ starts a gap $\equi$ $\gamma$ is either
  $\Pi_{3}$-reflecting, or is an admissible limit of extendibles. 
\end{corollary}

\section{A question of Rin on $\delta$-ITTMs}

It is obvous that if one has an ITTM with an uncountable number of cells on
the tape, {\eg}an $\omega_{1}$-ITTM, then one cannot address all of the cells
of the tape unless one allows ordinals parameters on the input tape of such
machines: there are only countably many programs. Taking the hint about the
uncountability of tape lengths causing unreachability through the scarcity of
programs, we give a direct argument that globally characterises the category
of tape lengths and their machines for which all cells can be reached.

We must turn then to considering computations in the class of functions
computable by a $\delta$-ITTM without parameters. We shall denote here the
class of such machines by $\delta$-ITTM$_{0}$. It is reasonable to require a
modest amount of closure on the length $\delta$ of the three tapes: we take
here closure under primitive recursive set functions (although weaker systems
would suffice). This ensures that $\delta$ is closed under G\"odel pairing.

\begin{theorem}
  \label{Th3.1}Every cell can be reached by a $\delta$-ITTM$_{0}$ computation
  if and only if during a run of a program, a wellorder of $\omega$ appears on
  (one of) its tape(s) at some stage, which has order type $\delta$.
\end{theorem}

Note: the requirement is not that the wellorder must be {\tmem{output}} by the
machine, (although we see that it can be) but only that it appears at some
point during the computation (in the description of Hamkins and Lewis
{\cite{HL}}, it is `accidental'). There is thus some computation on a
$\delta$-ITTM$_{0}$ which sees, even if temporarily, that `$\delta$ is
countable'. For $\delta$-ITTM$_{0}$ machines, Rin defines $\lambda_{0} (
\delta ) , \gamma_{0} ( \delta ) , \zeta_{0} ( \delta )$ {\etc}by analogy
with the ordinals for ITTMs, again see {\cite{HL}}, which we should call here
$\lambda_{0} ( \omega ) , \gamma_{0} ( \omega ) , \zeta_{0} ( \omega )$ ... \
(Note that there is no difference between functions $\omega$-ITTM$_{0}$
computable and $\omega$-ITTM computable, since any cell number `$n$' that
could be used as a parameter, is readily definable.)

\begin{definition}

  $\lambda_{0} ( \delta ) = \sup \{   || y ||: y \sset
  \omega \nocomma ,y \in \tmop{WO} \wedge y $ is the
  output of a halting $\delta$-ITTM$_{0}$ program$\nobracket \}$.
  
  $\Sigma_{0} ( \delta ) = \sup \{  || y || : y \sset
  \omega \nocomma ,y \in \tmop{WO} \wedge y $ appears
  at some stage on some tape of a $\delta$-ITTM$_{0}$ program$\nobracket \}$.
\end{definition}

We shall have also to consider the standard class of machines with length of
tape $\delta$: we'll write these without the zero subscript to indicate that
parameters are allowed. For example the definition of accidentally writable
wellorders appearing on such computations reads:

\begin{definition}
  \label{Def3.3}$\Sigma ( \delta ) = \sup \{ || y || :
  y \sset \omega \nocomma ,y \in \tmop{WO} \wedge y $
  appears at some stage on some tape $\delta$-ITTM program$\nobracket \}$
\end{definition}

We put to use primitive recursive closure of $\delta$, as in the next example to divide up
our scratch tape into $\delta$ many slices $\pa{S_{\alpha}}_{\alpha < \delta}$
with each $S_{\alpha}$ of order type $\delta$. On this array (call it $S$, and
we'll call the routine we describe next the $S$-{\tmem{routine}}) we write in
$\delta$ steps, {\via}a $\delta$-ITTM$_{0}$ computation, in the $\beta$'th
slice $S_{\beta}$, a code for $\beta$ (either as a real, or as a subset of
another ordinal, or even just a $\beta$-string of $1$'s followed by $0$'s, for
definiteness, let us say the latter). We arrange this in such a way that a
head never moves left whilst on a limit cell, but it first goes back to
$C_{0}$ when this routine is completed.

We shall use the $S$-routine as a prefix to a run of the
$\delta$-ITTM$_{0}$-universal program $U_{0}$, in which one subcomputation
actually simulates a run of the full $\delta$-ITTM$_{}$-universal program
$U$\tmtextbf{}: we may imagine, running $U_{0}$, that we simulate all
$\delta$-ITTM programs on $\omega \times \delta$ many slices of a scratch
tape initiating in turn all $U ( e, \tau )$ for some $e \in \omega ,$ and
increasing in turn some parameter $\tau < \delta$ which we now have access to
{\via}$S_{\tau}$. Notice then that any real, (indeed, we remark for below, any
subset of $\delta$) that appears on some \ $U ( e, \tau )$ tape also appears
coded on some subprogram of $U_{0}$ on the $  \delta$-ITTM$_{0}$ hardware.
(And trivially also the converse.) Hence

(1) $\Sigma_{0} ( \delta ) = \Sigma ( \delta ) \geq \delta$.

{\rem} Thus the difference between $\delta$-ITTM$_{0}$ and
$\delta$-ITTM-computations, taken as a class, is not that they reach different
classes of sets or functions during computations, it is the ability to halt
with, or output those results that creates a difference: {\eg}, an $f: \delta
\imp \delta$ may appear, but if for $\delta$-ITTM$_{0}$'s the cell $\alpha$
cannot be reached, we may not, in general have access to $f ( \alpha )$.

{\pf}of Theorem \ref{Th3.1}. $( \Leftarrow )$ Assume $\delta$-ITTM$_{0}$-machines cannot
reach some least cell $C_{\alpha_{0}}$ for some $\alpha_{0} < \delta$.

(2) $\lambda ( \delta ) \leq \alpha_{0} < \delta$.

{\pf}As Rin remarks (his Prop. 2.10): if a real coding $\beta$ is writable
then every cell $C_{\beta'}$ for $\beta' \leq \beta$ is addressable. {\qed}
(2)

(3){\tmem{ No function collapsing $\delta$ to $\omega$  appears on any
$\delta$-ITTM$_{0}$-computation.}}

{\pf} Using the slices $S_{\alpha}$ from the array $S$, or otherwise, we can
in a $\delta$-ITTM$_{0}$-computation build up a wellorder on $\delta \times
\delta$ of order type $\delta$. Thence wellorders of type $\delta + \delta$,
$\delta \cdot \omega , \ldots$ \ So order types of length $>  \delta$ can
appear on $\delta$-ITTM$_{0}$-computation tapes. So suppose for a
contradiction during a run of the universal $\delta$-ITTM$_{0}$ computation
some real $\bar{x} \in \tmop{WO}$ appears with $| \bar{x} | = \delta
\nobracket \nobracket$. (Thus showing $\Sigma ( \delta ) > | \bar{x} | =
\delta \nobracket \nobracket$.) We arrange our run of this universal
computation, so that whenever a real $x' \in \tmop{WO}$ appears, we split off
to a subroutine that runs, for example, a copy of the $S$-routine indexing its
stages along the wellorder $x'$ until either $x'$ is exhausted or $S$ is
completed. If the former, it returns to the main computation; if the latter it
HALT's with output $x'$. However then $\lambda ( \delta ) \geq | x' | \geq
\delta \nobracket \nobracket$, contrary to (2).{\qed}(3)

$( \Rightarrow )$ For the converse assume $\delta$ is such that every cell can
be reached by a $\delta$-ITTM$_{0}$-program. During a run of $U$ we build a
wellorder $E$ on a scratch tape: let $e ( \alpha )$ be that $e \in \omega$ for
which $P_{e}$ halts first with its head on $C_{\alpha}$. As we run $U$ we set

$n E m  \Equi \ex \alpha < \alpha'   ( n=e ( \alpha ) \wedge m=e ( \alpha' )
) .$

{\nod}Then $E \in \tmop{WO}$, and demonstrates that $\delta \leq \lambda (
\delta )$: it will be completed and is even writable, once all the relevant
$P_{e}$'s have halted. As $E$ collapses $\delta$ this direction and the
theorem are complete.

{\qed} Thm \ref{Th3.1}\\

{\rem}As {\cite{HL}} shows, if $\alpha$ appears (coded)
on a tape, there is another computation with a code for $L_{\alpha}$ appearing
on a tape. \ Then, by absoluteness of computations to the $L$-hierarchy, we
shall have:

\begin{corollary}
  \label{Cor3.4}{\tmstrong{}} \tmverbatim{}Every cell can be addressed in the
  $\delta$-ITTM$_{0}$-class if and only if $L_{\Sigma ( \delta )}
  \models$``$\delta  $ is countable''.
\end{corollary}

We close this section with a remark {\em \`a propos} the Introduction: with the formalism of constants in the language we can now formulate the theory of quasi-inductive definitions over the structure $(\delta, <)$ or indeed over a $\delta$-ITTM considered as  a relational structure generalising Moschovakis in \cite{M}.  Without such we should be restricted to structures satisfying the conditions of the last Corollary, which would be unnatural as well as undesirable.

{\small{}}\section{Infinite Time Blum-Shub-Smale machines}

This section seeks to show that the strength of IBSSM's equipped with a Liminf
(in the Euclidean metric) rule for the register contents, rather than the
original continuity requirement used in {\cite{KoSe12}} allows,
unsurprisingly, for much greater computational power. It is perhaps more surprising
that ITTMs can be simulated on such a machine, so they have the latter's power at least. It is easy then to see that the power is exactly that of
ITTMs: since ITTMs can construct $L$ and as IBSSM's computations, on
rational input, are absolute to $L$, the latter can compute no more than the
former.

We assume an IBSSM program $B_{e}$ with registers $R_{0} ,R_{1 \nocomma \nocomma
,}   \ldots . \nocomma ,R_{n}$. $B$ has finitely many instructions $I_{0}
\nocomma , \nocomma \cdots \, \nocomma ,I_{k}$. We suppose that at time
$\alpha$ the instruction number in the enumeration of the flow diagram of $B_{e}$
about to be effected is $I ( \alpha \nobracket$); again at limit times
$\lambda$, we set $I ( \lambda )$ to be  $l= \liminf_{\beta \rightarrow
\lambda} I ( \beta )$ as in {\cite{KoSe12}}.

\begin{definition}
  Let $R_{i} ( \alpha ) \in \re$ be the real contents of register $R_{i}$ at
  time $\alpha$. 
\end{definition}

We shall assume, wlog, that a real is expressed by its decimal expansion. At
limit stages of time $\lambda$ $R_{i} ( \lambda ) \dfs \tmop{Liminf}_{\alpha
\rightarrow \lambda} R_{i} ( \alpha ) .$ (If this Liminf becomes infinite we
make some other decision, for sake of definiteness, allow the machine to
crash.) Assume we start with 0.0's (or rationals) in all registers.
Note that for any limit $\lambda$, $i\leq
n$, $u=\pa { R_{i}(\gamma)\mid
\gamma < \alpha}$ is 
$\Sigma_{1}^{L_{\lambda}}(\{\alpha\})$ for any $\alpha< \lambda$.
 Let $(
\zeta , \Sigma )$ be as defined for ITTMs.

\begin{lemma}
  $R_{i} ( \zeta ) =R_{i} ( \Sigma )$, and $B_{e}$ has either halted or, as for
  {\tmem{ITTM}}s, starts looping by at most stage \nolinebreak $\zeta$.
\end{lemma}

{\pf} 
We let $Q^{n}$ be the collection of rationals whose decimal expansion is of the form  $0.r_{1}r_{2}\cdot
r_{n}00000$.
Let $s^{i}_{n} ( \alpha )\in Q^{n}$ be the rational approximation whose expansion
is the first $n$ digits of $R_{i} ( \alpha )$.
 We'll simplify the discussion by assuming that a) we
discuss only $R_{0}$ throughout and \ b) $ R_{0} ( \alpha ) \in [ 0,1 )$.
Towards this end we first define for limit ordinals $\lambda$:
$$I^{+}( t,n, \lambda ) \equi 
 t \in Q^{n} \wedge\all \alpha <
\lambda \ex \beta \in ( \alpha , \lambda ) (t > R_{0} ( \beta ) \wedge t- R_{0}(\beta) \leq 10^{-n} ).$$
If $R_{0}(\lambda)\da$ (that is, is less than $\infty)$) then $\ex t \all n I^{+}(t,n,\lambda)\neq
\emp$. Let $s_{n,\lambda}^{+} = \min \{t\mid  I^{+}(t,n,\lambda)\}$. Then $s_{n,\lambda}^{+} = y \equi y\in I^{+}(t,n,\lambda)\}\wedge y -10^{-n}\notin I^{+}(t,n,\lambda)\}$. Hence $s_{n,\lambda}^{+}$ is $\Pi_{2}\wedge\Sigma_{2}$ definable.\\

\nod {\em Claim}\quad  $s_{n,\lambda}^{+} - 10^{-n}\leq s^{0}_{n}, R_{0}(\lambda)\leq s_{n,\lambda}^{+}$.

\pf of {\em Claim}. To avoid too much clutter, we set $s= s^{0}_{n} ( \lambda)$ and $s^{+}= s_{n,\lambda}^{+}$. We divide into two cases depending on how the liminf value at $R_{0}(\lambda)$ can be approached.\\

\nod {\em Case 1} $R_{0}(\lambda)$ {\em is a simple limit from below of a subsequence $\pa{ R(\alpha_{\ii})\mid \ii < \tau}$ for some $\tau$ with $\la \alpha_{\ii}\ra_{\ii<\tau} $ cofinal in } $\lambda$.

{\em Case 1a)} $R_{0}(\lambda)=s$.   Then one can check that our definitions imply that $s^{+}=s$ and $s^{+}=R_{0}(\lambda)>s^{+}-10^{-n}$.

{\em Case 1b)} $R_{0}(\lambda)>s$.  Then  $s^{+}> R_{0}(\lambda)>s^{+}-10^{-n}= s$.\\
Hence in both subcases we have $s^{+}\geq R_{0}(\lambda)>s^{+}-10^{-n}$ with the equality holding if and only if $R_{0}(\lambda)=s$.\\

\nod {\em Case 2} {\em The sequence $\pa{ R_{0}(\alpha)}$ for $\alpha < \lambda$ approaches $R(\lambda)$ purely as a liminf from above (\ie from some point on, no $R_{0}(\alpha)$ is less than $R_{0}(\lambda)$).}

Then, similar to the first case,  $s^{+}> R_{0}(\lambda), s \geq s^{+}-10^{-n}$, with $R_{0}(\lambda) = s^{+}-10^{-n}$ holding if and only if $R_{0}(\lambda)=s$. \qed ({\em Claim})\\

\nod From the {\em Claim} it follows that $$ (\ast) \quad \all \lambda, \lambda' [\all n(s_{n,\lambda}^{+}=s_{n,\lambda'}^{+}) \imp R_{0}(\lambda)=R_{0}(\lambda')].$$

{\nod}If we imagine the computation proceeding by a recursion within the $L$
hierarchy then we shall have that $R _{0}( \beta )$ is a real definable over
$L_{\beta}$, by using, for example, the set of  rational approximations $\{
s^{+}_{n, \beta } \}_{n< \omega}$. 
 As
$L_{\zeta} \prec_{\Sigma_{2}} L_{\Sigma}$ we shall have for any $n,t$ that $I
( t,n, \zeta ) \equi I ( t,n, \Sigma )$ and thus $s_{n,\zeta}^{+}=s_{n,\Sigma}^{+}$. But then $R_{0} ( \zeta ) =R_{0} (
\Sigma )$ by $(\ast)$ above. Identical arguments hold for the other registers $R_{i}$. Also then by another straightforward application of
$\Sigma_{2}$-reflection 
we have that  $R_{i} ( \gamma ) \geq R_{i} ( \zeta )$ for all $\gamma \in [
\zeta , \Sigma ]$.


These facts imply that no Liminf smaller than $R_{i} (
\zeta )$ can ever appear in $R_{i}$ at a later part of the computation, and the
computation is looping. (This can also be seen by observing the $B_{e}$
computation coded inside the ITTM-theory machine - see below.) \ {\qed} \\

The harder part is the lower bound. A natural attempt is to somehow define an
ITBSS program that continuously runs producing (reals coding) the theories of
successive levels of the $L$-hierarchy much as in {\cite{FrWe08}}. In, {\eg}\!\!, {\cite{We2015}} at Theorem 11 we stated how one might produce successive
Turing jumps for any $\alpha < \omega^{\omega}$ on the continuity-IBSSM's.
That however required one to know in advance how long a hierarchy was
involved. In {\cite{FrWe08}} the ITTM runs continuously irrespective of any
fixed length, and indeed loops forever, deliberately so. It would seem
difficult (but presumably not impossible in view of the sequel) to formulate
this directly as an IBSSM. Instead we show how to simulate an ITTM computation
$P_{e} ( 0 )$ directly on an IBSSM as some $B_{\widebar{e}} ( 0.0 )$.

The obvious first, and really only, hurdle if we don't know how long the
computation is intended to run, is how to keep a record of the infinite tape
contents. Indeed, since the rest of the ITTM's actions are simple effects from
its state-transition table, which are easily simulable as IBSSM actions,
keeping the ITTM cell tape's contents is the main difficulty. \ It is easy to
see what goes wrong with a too-simple approach: suppose one has a real with a
decimal expansion 0.1101... with a 1 in the $k$'th position iff $C_{k}$
contains a $1$; suppose before stage $\omega$ cells $C_{0}$ and $C_{1}$ have
both changed value infinitely often, but so that at no finite stage are they
both 0, then by the ITTM-Liminf rule $C_{0} ( \omega ) =C_{1} ( \omega ) =0$.
But the IBSSM-Liminf rule applied to the sequence of reals as decimals of
the form $0.10 \tmop{xxx}$ $0.01 \tmop{xxxx}$ always with a 1 in the first or
second position, yields a liminf of $0.01 \tmop{xxx}$ rather than the
intended $0.00 \tmop{xxxxx}$.

So instead we use not just one, but an infinite number of places in the real's
expansion to code the correct cell contents of {\egc} $C_{0}$.

Let $\pi : \omega \times \omega\,\, \lr \,\,\omega$ be a recursive bijection. Let
$H_{i} = \left\{ n \mid \ex k \pi \left( \pa{i,k} \right) =n \right\}$. Let $h_{i} :
\omega \twoheadrightarrow H_{i}$ be the, also recursive, strictly increasing
enumeration of $H_{i}$. We use $H_{i}$ to record information about $C_{i}$.
Suppose w.l.o.g. we initiate the ITTM-computation $P_{e} ( 0 )$ on zero input
with zeros in every cell.   In the sequel, $r(\alpha)$ denotes the real that is in $R_{0}$ at time $\alpha$.  We start the IBSSM-simulation $B_{e} ( 0.0000 \ldots
)$ by placing a real $r ( 1 ) \in [ 0,1 )$ in the register $R_{0}$ whose decimal
expansion has 0 everywhere except for the $h_{i} ( 0 )$ position for each $i<
\omega$, where it has instead a 1. ($r ( 1 )$ thus has at time $t=1$,
infinitely many 1's in its expansion.)

The idea is that at stage $\alpha$ \ $C_{i} ( \alpha )$ equals a zero (or one)
iff there is an even (resp. odd) number $k$ with the $h_{i} ( k )$'th position
a 1 in the decimal expansion of $r( \alpha )$. For this to work at every
moment in time $\alpha$, each of the $H_{i}$ sets will have at most one 1 at
its corresponding point in $r ( \alpha \nobracket$). Otherwise the positions
in $H_{i }$ are all 0.

{\bu} Note first that given a real $r$ in $R_{0}$, there are IBSSM-operations
that first check whether $r$ is conformable to a code  for the whole tape
$\vec{C_{i}} \nocomma$; if so it can then read off whether $C_{i}$ is zero (or
one), as coded into $r$: this involves computing whether $r$ has in its
expansion a 1 at position $h_{i} ( 2k )$ (or at $h_{i} ( 2k+1 )$). These
latter checks merely involve calculations of the recursive function $\pi$,
multiplications and divisions by powers of 10 {\etc}\!, {\etc}and can be
programmed into a subroutine of the IBSSM $B_{e}$'s flow diagram. (Here and
later we leave all this as an exercise for the reader).

{\bu} Secondly, suppose at time $\alpha$ we have $C_{i} ( \alpha ) =1 \neq
C_{i} ( \alpha +1 )$. There is thus a change in the value of $C_{i}$ initiated
by the ITTM $P_{e}$. We effect the following change on $r( \alpha +1 )$: we move
the current $1$ in its expansion,  at some position of the form
$h_{i} ( 2k+1 )$ to position $h_{i} ( 2k+2 )$, and the $h_{i} ( 2k+1 )$'th
position in what will be $r ( \alpha +2 )$ reverts back to 0.  Again these
changes to $r ( \alpha +1 )$ in $R_{0}$ can be effected by simple IBSSM
programmable arithmetic. (For $C_{i} ( \alpha ) =0 \neq C_{i} ( \alpha +1 )$
we make similar changes {\mm}.)

The reader may wish to {\egc} consider what happens to just cells $C_{0}$ and
$C_{1}$ at some limit stage $\lambda$ in $P_{e}$. Suppose $C_{0}$ changes
value cofinally often below $\lambda$ but all other cells are fixed from some
point $\alpha$ onwards below $\lambda$: each time there is a change a 1 gets
pushed down a further stage in the $H_{0}$-numbered positions to the right,
thus decreasing  the real in $R_{0}$ at least in terms of the values at the positions in
$H_{0}$. After infinitely many such moves to the right, when we consider $r (
\lambda )  =  \tmop{Liminf}_{\beta \rightarrow \lambda} r ( \beta )$, this
will be simply the real $r ( \alpha )$ but with all of the $H_{0}$
positions now set to zero. This corresponds to having in the ITTM-computation the Liminf
operation setting $C_{0} ( \lambda )$ to 0. To mirror this, at stage $\lambda
+1 $ we simply modify the real $r(\lambda)$ in $R_{0}$ by just resetting the $h_{0} ( 0 )$'th decimal position to a 1 (just as we initially did for $r(1)$) ready
for any further changes at $C_{0} ( \gamma \nobracket$) for $\gamma >
\lambda$, that may occur (and leave all else alone).

Suppose $C_{1}$ alters infinitely often below $\lambda$ as well. Then all the
$H_{1}$ positions in the Liminf real $r ( \lambda )$ have become 0 also. As
before $r ( \lambda +1 )$ is calculated so that the $h_{1} ( 0 )$ position
contains a 1.

However now note the relative independence of the IBSSM Liminf operation:
making these infinitely many changes to $r$ caused by the changes of $C_{0}$
followed afterwards by infinitely many changes to $r$, caused by those of
$C_{1}$ sequentially, (all other cells being kept constant) results in the
same Liminf real as if the changes to $C_{0}$ and $C_{1}$ had happened
concurrently. If before a limit $\lambda$ only finitely many changes have
occurred in $C_{j}$ since there was last a 1 at position $h_{j} ( 0 )$, then
$C_{j} ( \lambda ) =C_{j} ( \beta )$ where the last change occurred at time
$\beta$. This $1$ is then preserved into the Liminf $r( \lambda )$.

{\bu} Thus, to sum up, at any limit stage $\lambda$ after $r ( \lambda ) =
\tmop{Liminf}_{\alpha \rightarrow \lambda} r ( \alpha )$ has been set, $B_{e}$
then inspects each of the $H_{i}$ sets in the expansion of $r ( \lambda )$ in
turn, and if $\tmop{ran} ( h_{i} )$ has no $1$ in its range, then it
calculates $r ( \lambda +1 )$ so that $h_{i} ( 0 ) =1$; it does this for each
$i< \omega$ in turn (possibly using further registers as scratch area). As
before these are IBSSM-calculable arithmetic operations.

{\bu} A record of the head position of $P_{e}$ can be recorded in a register
$R_{2}$ of $B_{e}$: if this R/W head of $P_{e}$ is over the triple of cells
$\nobracket C_{3k} \nocomma ,C_{3k+1} ,C_{3k+2} \nobracket  $ then set $R_{2}$
to $\frac{1}{k}$ (for $k>0$, for $k=0$ let $R_{2}$ be 0). Then if our model of
ITTM computation replaces the R/W head to the beginning of its tape at each
limit stage, we can set further flag registers to alert us to limit stages in
the $P_{e}$ computation, in order to register this fact in $R_{2}$.

This finishes (modulo some other programming details of no great interest)
the simulation of $P_{e} ( 0 )$ as some $B_{e} ( 0 )$. {\qed}

\begin{corollary}
  The {\tmem{ITTM-}} and {\tmem{IBSSM-}}models are ``bi-simulable''.
\end{corollary}

We have just prescribed a computation $B_{\bar{e}} ( 0.00 \ldots )  $ from
$P_{e} ( 0 )$. But any given $B_{f} ( 0 )$ can be run recursively inside the
$L$-hierarchy, with a snapshot of its $\alpha$'th stage definable over
$L_{\alpha}$. \ However the theory machine of {\cite{FrWe08}} constructs
levels of the hierarchy, and their theories, for $\alpha < \Sigma$. It can
thus decode the action of $B_{f}$ and find the sequence of snapshots of
$B_{f}$. This provides a way of emulating $B_{f}$ in some ITTM computation
$P_{f'}$. {\qed}

\begin{corollary}
  \label{Cor20}There is a universal {\tmem{IBSSM}}.
\end{corollary}

{\pf} Every IBSSM program can be simulated on an ITTM. Because there is a
universal ITTM program $U$ say, and in turn this program can be simulated
on an IBSSM as discussed above, we can design a universal IBSSM that
decodes from the simulation of $U$ on an IBSSM as described above, the action
of any desired $B_{e} ( k )$ for any IBSSM-program index $e$. {\qed}\\

Now we have such a result, we can take over into the Liminf-IBSSM many of the
results and ideas from ITTM-theory. (For example we now know the sup of the
IBSSM clockable or writable or eventually writable ordinals.)

We then have:

\begin{corollary}
  \label{C4.5} For {\tmem{IBSSM}} with the $\liminf$rule, then the set $Z$ of
  reals on which an {\tmem{IBSSM}} computation (on rational input) is convergent, is
  precisely the set $Z $of those reals in $L_{\zeta}$.
\end{corollary}

A {\tmem{convergent computation}} (on rational input) is one in which a
designated register, {\eg}$R_{0}$ has a final settled value from some point
onwards.

It would be of interest to find questions or properties more in the
geometric/Euclidean flavour of the BSS-idea.

Let $\mathbbm{T}$ be the $n$-dimensional torus obtained as $[ 0,1 )^{n}$ (with
$1$ identified with $0$). Let $f$ be a continuous or effective function,
$f:\mathbbm{T} \rightarrow \mathbbm{T}$. Let $f^{\alpha}$ be the $\alpha$'th
iterate of $f.$ Here $f^{\alpha +1} ( x )  = f ( f^{\alpha} ( x ) )$, whilst
$f^{\lambda} ( x ) =  \liminf_{\beta \rightarrow \lambda} f^{\beta} ( x )$ if
$\tmop{Lim} ( \lambda )$. By the above an IBSSM can compute $f^{\alpha} ( x )$
for any $x \in \re$. Let $\mathcal{O}= \left\{ p \in \mathbbm{T} \mid   \ex
\alpha   ( f^{\alpha} ( p ) =0 ) \right\}$ ($\mathcal{O}$ is the `origin set',
the set of points that at some iterate of $f$ get sent to the origin $0$). One
may show the following.

\begin{theorem}
  $\Pi^{1}_{3}$-$\tmop{CA}_{0}$ proves the existence of the set $\mathcal{O}$.
  $\Pi^{1}_{2}$-$\tmop{CA}_{0}$ \ is insufficient.
\end{theorem}

\begin{theorem}
  Let $\tau_{f} = \sup \left\{ \alpha \mid \ex p \in \mathbbm{T} \right.
  \nobracket ( f^{\alpha} ( p ) =0 ) \mbox{ but } \all \beta < \alpha (f^{\beta} ( p ) \neq 0 ) \}$, for recursive $f$. \ Then 
  
  (i)
  $\tau_{f} \leq \lambda$.
  
  (ii) There exists a recursive $f$ that attains this bound. 
\end{theorem}

This is just a sample of possible properties from topological dynamics that are open to investigate when defined with transfinite actions.

\small
\bibliographystyle{plain}
\bibliography{settheory10t}

\begin{thebibliography}{10}

\bibitem{PH2}
P.~Aczel and W.~Richter.
\newblock Inductive definitions and analogies of large cardinals.
\newblock In W.~Hodges, editor, {\em Conference in Mathematical Logic London
  70}, volume 255 of {\em Lecture Notes in Mathematics}, pages 1--9. Springer,
  1971.

\bibitem{Bar}
K.J. Barwise.
\newblock {\em Admissible Sets and Structures}.
\newblock Perspectives in Mathematical Logic. Springer Verlag, 1975.

\bibitem{BeBuFr12}
A.~Beckmann, S.~Buss, and S-D. Friedman.
\newblock Safe recursive set functions.
\newblock {\em J. Symbolic Logic}, 80(3):335--369, 2015.

\bibitem{BL80}
A.~Beller and A.~Litman.
\newblock A strengthening of {Jensen's} {${\Box}$} principles.
\newblock {\em Journal of Symbolic Logic}, 45(2):251--264, 1980.

\bibitem{Burd84}
B.~Burd.
\newblock The speed of an infinite computation.
\newblock Master's thesis, Rutgers, May 1984.

\bibitem{B}
J.P. Burgess.
\newblock The truth is never simple.
\newblock {\em Journal of Symbolic Logic}, 51(3):663--681, 1986.

\bibitem{Carl2019}
M.~Carl.
\newblock {\em Ordinal {Computability}: An {Introduction} to {Infinitary}
  {Machines}}.
\newblock De Gruyter Series in Logic and Its Applications. De Gruyter, 2019.

\bibitem{Carl20}
M.~Carl.
\newblock Clockability for {Ordinal} {Turing } {Machines}.
\newblock Lecture Notes in Computer Science. Springer, 2020.

\bibitem{CaRiSc18}
M.~Carl, B.~Rin, and P.~Schlicht.
\newblock Writability and reachability for $\alpha$-tape infinite time {Turing}
  machines.
\newblock {\em ArXive}, 1612:02982v1, February 2018.

\bibitem{De}
K.~Devlin.
\newblock {\em Constructibility}.
\newblock Perspectives in Mathematical Logic. Springer Verlag, Berlin,
  Heidelberg, 1984.

\bibitem{D}
A.~J. Dodd.
\newblock {\em The {Core Model}}, volume~61 of {\em London Mathematical Society
  Lecture Notes in Mathematics}.
\newblock Cambridge University Press, Cambridge, 1982.

\bibitem{FrWe08}
S-D Friedman and P.D. Welch.
\newblock Two observations regarding {Infinite} {Time} {Turing} machines.
\newblock In I.~Dimitriou, editor, {\em Bonn Interational Workshop on Ordinal
  Computability}, pages 44--48, Bonn, 2008. Hausdorff Centre for Mathematics,
  University of Bonn,
  http://www.math.uni-bonn.de/ag/logik/events/biwoc/report.pdf.

\bibitem{Ga74}
R.O. Gandy.
\newblock Set theoretic principles for elementary syntax.
\newblock In T.~Jech, editor, {\em Axiomatic Set Theory}, volume 13 II of {\em
  Proceedings of Symposia in Pure Mathematics}, pages 103--126, Providence,
  Rhode Island, 1974. American Mathematical Society.

\bibitem{GB}
A.~Gupta and N.~Belnap.
\newblock {\em The revision theory of truth}.
\newblock M.I.T. Press, Cambridge, 1993.

\bibitem{HL}
J.D. Hamkins and A.~Lewis.
\newblock Infinite time {Turing} machines.
\newblock {\em Journal of Symbolic Logic}, 65(2):567--604, 2000.

\bibitem{HaSe}
J.D. Hamkins and D.~Seabold.
\newblock Infinite time {Turing} machines with only one tape.
\newblock {\em Mathematical Logic Quarterly}, 47(2):271--287, 2001.

\bibitem{JeKa71}
R.B. Jensen and C.~Karp.
\newblock Primitive recursive set functions.
\newblock In D.~Scott, editor, {\em Axiomatic Set Theory}, volume 13 I of {\em
  Proceedings of Symposia in Pure Mathematics}, pages 143--167, Providence,
  Rhode Island, 1971. American Mathematical Society.

\bibitem{Kl62b}
S.C. Kleene.
\newblock Turing-machine computable functionals of finite type {I}.
\newblock In {\em Proceedings 1960 Conference on Logic, Methodology and
  Philosophy of Science}, pages 38--45. Stanford University Press, 1962.

\bibitem{Kl62a}
S.C. Kleene.
\newblock Turing-machine computable functionals of finite type {II}.
\newblock {\em Proceedings of the London Mathematical Society}, 12:245--258,
  1962.

\bibitem{K05}
P.~Koepke.
\newblock Turing computation on ordinals.
\newblock {\em Bulletin of Symbolic Logic}, 11:377--397, 2005.

\bibitem{KoSe12}
P.~Koepke and B.~Seyfferth.
\newblock Towards a theory of infinite time {Blum-Shub-Smale} machines.
\newblock In S.~Cooper, A.~Dawar, and B.~L\"owe, editors, {\em How the World
  Computes}, volume 7318 of {\em Lecture Notes in Computer Science}, pages
  405--415. Spinger, 2012.

\bibitem{L}
B.~L{\"{o}}we.
\newblock Revision sequences and computers with an infinite amount of time.
\newblock {\em Journal of Logic and Computation}, 11:25--40, 2001.

\bibitem{Mach61}
M.~Machover.
\newblock The theory of transfinite recursion.
\newblock {\em Bulletin of the American Mathematical Society}, 67:575--578,
  1961.

\bibitem{M}
Y.N. Moschovakis.
\newblock {\em Elementary Induction on Abstract structures}, volume~77 of {\em
  Studies in Logic series}.
\newblock North-Holland, Amsterdam, 1974.

\bibitem{Ri79}
T.L. Richardson.
\newblock {\em A {Silver} Machine approach to the Constructible Universe}.
\newblock PhD thesis, University of California, Berkeley., 1979.

\bibitem{Rin}
B.~Rin.
\newblock The computational strengths of $\alpha$-tape infinite time {Turing}
  machines.
\newblock {\em Annals of Pure and Applied Logic}, 165(9):1501--1511, 2014.

\bibitem{Rog}
H.~Rogers.
\newblock {\em Recursive Function Theory}.
\newblock Higher Mathematics. McGraw, 1967.

\bibitem{Si99}
S.~Simpson.
\newblock {\em Subsystems of second order arithmetic}.
\newblock Perspectives in Mathematical Logic. Springer, January 1999.

\bibitem{Tak60}
G.~Takeuti.
\newblock On the recursive functions of ordinal numbers.
\newblock {\em Journal of the Mathematical Society of Japan}, 12:119--128,
  1960.

\bibitem{W}
P.D. Welch.
\newblock Eventually {Infinite} {Time} {Turing} degrees: infinite time
  decidable reals.
\newblock {\em Journal of Symbolic Logic}, 65(3):1193--1203, 2000.

\bibitem{W1}
P.D. Welch.
\newblock On revision operators.
\newblock {\em Journal of Symbolic Logic}, 68(3):689--711, 2003.

\bibitem{W09}
P.D. Welch.
\newblock Characteristics of discrete transfinite {Turing} machine models:
  halting times, stabilization times, and normal form theorems.
\newblock {\em Theoretical Computer Science}, 410:426--442, January 2009.

\bibitem{W2010}
P.D. Welch.
\newblock Discrete transfinite computation models.
\newblock In S.~B. Cooper and A.~Sorbi, editors, {\em Computability in
  {Context}: {Computation} in the real world}, pages 375--414. Imperial College
  Press/World Scientific, 2010.

\bibitem{We2015}
P.D. Welch.
\newblock Discrete transfinite computation.
\newblock In G.~Sommaruga and T.~Strahm, editors, {\em Turing's revolution :
  the impact of his ideas about computability}. Birkh\"auser Verlag, Basel,
  2015.

\end{thebibliography}

\end{document}